# ON THE SCALING OF THE CHEMICAL DISTANCE IN LONG-RANGE PERCOLATION MODELS

By Marek Biskup

*University of California, Los Angeles*

We consider the (unoriented) long-range percolation on $\mathbb{Z}^d$ in dimensions $d \geq 1$, where distinct sites $x, y \in \mathbb{Z}^d$ get connected with probability $p_{xy} \in [0, 1]$. Assuming $p_{xy} = |x - y|^{-s + o(1)}$ as $|x - y| \to \infty$, where $s > 0$ and $|\cdot|$ is a norm distance on $\mathbb{Z}^d$, and supposing that the resulting random graph contains an infinite connected component $\mathscr{C}_\infty$, we let $D(x, y)$ be the graph distance between $x$ and $y$ measured on $\mathscr{C}_\infty$. Our main result is that, for $s \in (d, 2d)$,

$$D(x, y) = (\log|x - y|)^{\Delta + o(1)}, \qquad x, y \in \mathscr{C}_\infty, \ |x - y| \to \infty,$$

where $\Delta^{-1}$ is the binary logarithm of $2d/s$ and $o(1)$ is a quantity tending to zero in probability as $|x - y| \to \infty$. Besides its interest for general percolation theory, this result sheds some light on a question that has recently surfaced in the context of "small-world" phenomena. As part of the proof we also establish tight bounds on the probability that the largest connected component in a finite box contains a positive fraction of all sites in the box.

## 1. Introduction.

### 1.1. *Motivation.*

Percolation is a simple but versatile model with applications ranging from the study of phase transitions in mathematical physics to opinion spreading in social sciences. The most well-understood questions of percolation theory are those concerning the appearance and uniqueness of the infinite component [10], uniqueness of the critical point [1, 16, 20], decay of connectivity functions [11, 12], and the scaling properties at the critical point in dimensions $d = 2$ [24, 25] and $d$ large enough [17, 18]. Less well understood remain natural questions about the qualitative structural and geometrical properties of the infinite connected component, especially below the upper critical dimension. In particular, this includes the tantalizing open problem concerning the absence of percolation at the percolation threshold.

Long-range versions of the percolation model have initially been introduced in order to study the effect of long-range interaction on the onset of phase transition in one-dimensional systems. On $\mathbb{Z}$, the most common setup is that, in addition to random nearest-neighbor connections with probability $p \in (0, 1)$, a bond between $x, y \in \mathbb{Z}$ is added with probability $1 - \exp\{-\beta|x-y|^{-s}\}$, where $\beta \in (0, \infty)$ and $s > 0$. In dimension 1, the interesting ranges of values of $s$ are $s < 1$, where the resulting graph is almost surely connected [23], $1 < s < 2$, where an infinite component appears once $p$ is large enough [22], and the critical case $s = 2$, where the infinite component appears "discontinuously" for some $p < 1$ sufficiently large if and only if $\beta > 1$ [2] and where the truncated connectivity function decays with a $\beta$-dependent exponent [19] for $\beta$ in the interval $(1, 2)$. The cases $s > 2$ are qualitatively very much like the nearest-neighbor case (in particular, there is no percolation for $p < 1$ and $\beta < \infty$). In dimensions $d > 1$, the insertion of long-range connections is not essential for the very existence of percolation—the main problem of interest there is to quantify the effect of such connections on the critical behavior.

In this paper we study the global scaling properties of the infinite component in long-range percolation models on $\mathbb{Z}^d$ for arbitrary $d$. We focus on the scaling of the graph distance (aka *chemical distance*) in the cases when the probability that a bond is occupied falls off with exponent $s \in (d, 2d)$. More precisely, we let distinct $x, y \in \mathbb{Z}^d$ be connected independently with probability $p_{xy}$ that has the asymptotics $p_{xy} = 1 - \exp\{-|x-y|^{-s+o(1)}\}$ as $|x-y| \to \infty$. Assuming that there is a unique infinite connected component $\mathscr{C}_\infty$ almost surely, we let $D(x, y)$ be the distance between the sites $x$ and $y$ measured on $\mathscr{C}_\infty$. Then we prove that $D(x, y)$ scales with the Euclidean distance $|x-y|$ as

$$(1.1) \qquad D(x, y) = (\log|x-y|)^{\Delta + o(1)}, \qquad x, y \in \mathscr{C}_\infty, \ |x-y| \to \infty,$$

where $\Delta = \Delta(s, d)$ is given by

$$(1.2) \qquad\qquad \Delta(s, d) = \frac{\log 2}{\log(2d/s)}.$$

This result should be contrasted with those of [[5]–[7, 14]] (see also [4]), where various (other) regimes of decay of long-range bond probabilities have been addressed. We refer to Section 1.3 for further discussion of related work and an account of the current state of knowledge about the asymptotic behaviors of $D(x, y)$.

The nonlinear scaling (1.1) is a manifestation of the fact that adding sparse (but dense-enough) long edges to a "Euclidean" graph may substantially alter the graph geometry and, in particular, its scaling properties. This is exactly what has recently brought long-range percolation renewed attention in the studies of the so-called "small-world" phenomena; see [26] for an initial work on these problems. This connection was the motivation of the work by Benjamini and Berger [6], who studied how the (graph) diameter of a finite ring of $N$ sites changes when long connections are added in. On the basis of a polylogarithmic *upper* bound, the authors of [6] conjectured (cf. Conjecture 3.2 in [6]) that in the regime when $s \in (d, 2d)$, the diameter scales as $(\log N)^\gamma$, where $\gamma = \gamma(s) > 1$. The present paper provides a polylogarithmic *lower* bound in this conjecture. However, at present it is not clear whether the exponent for the diameter growth matches that for the *typical* distance between two remote points. We refer to Section 1.4 for further discussion of "small-world" phenomena.

The remainder of this paper is organized as follows. In Section 1.2 we define precisely the long-range percolation model and state our main theorem (Theorem 1.1). In Section 1.3 we proceed by summarizing the previous results concerning the behavior of $D(x, y)$—and graph diameter—for various regimes of $s$. In Section 1.4 we discuss the relation to "small-world" phenomena. Section 2 is devoted to a heuristic explanation of the proof of Theorem 1.1. The proof requires some preparations, particularly an estimate on the size of the largest connected component in large but finite boxes. This is the content of Theorem 3.2 in Section 3. The actual proof of our main result comes in Sections 4.1 (upper bound) and 4.2 (lower bound).

1.2. *The model and main result.* Consider the $d$-dimensional hypercubic lattice $\mathbb{Z}^d$ and let $(x, y) \mapsto |x - y|$ denote a norm distance on $\mathbb{Z}^d$. For definiteness, we can take $|\cdot|$ to be the usual $\ell^2$-norm; however, any other equivalent norm will do. Let $q: \mathbb{Z}^d \to [0, \infty)$ be a function satisfying

$$(1.3) \qquad \lim_{|x| \to \infty} \frac{\log q(x)}{\log |x|} = -s,$$

where $s \geq 0$. (Here we set $\log 0 = -\infty$.) For each (unordered) pair of *distinct* sites $x, y \in \mathbb{Z}^d$, we introduce an independent random variable $\omega_{xy} \in \{0, 1\}$ with probability distribution given by $\mathbb{P}(\omega_{xy} = 1) = p_{xy}$, where

$$(1.4) \qquad p_{xy} = 1 - e^{-q(x-y)}.$$

Note that $p_{xy} = p_{x-y,0}$ so the distribution of $(\omega_{xy})$ is translation invariant.

Let $\mathscr{G}$ be the random graph with vertices on $\mathbb{Z}^d$ and a bond between any pair of distinct sites $x$ and $y$, where $\omega_{xy} = 1$. Given a realization of $\mathscr{G}$, let us call $\pi = (z_0, z_1, \ldots, z_n)$ a *path*, provided $z_i$ are all *distinct* sites in $\mathbb{Z}^d$ and $\omega_{z_{i-1}z_i} = 1$ for each $i \in \{1, 2, \ldots, n\}$. Define the length $|\pi|$ of $\pi$ to be the number of bonds constituting $\pi$ (i.e., the number $n$ above). Using $\Pi(x, y)$ to denote the (random) set of all paths $\pi$ with $z_0 = x$ and $z_{|\pi|} = y$, we let

$$(1.5) \qquad D(x, y) = \inf\{|\pi| : \pi \in \Pi(x, y)\}, \qquad x, y \in \mathbb{Z}^d.$$

[In particular, we have $D(x, y) = \infty$ if $\Pi(x, y) = \varnothing$.] The random variable $D(x, y)$ is the *chemical distance* between $x$ and $y$, that is, the distance measured on the graph $\mathscr{G}$.

Throughout the rest of the paper, it will be assumed that the random graph $\mathscr{G}$ contains an infinite connected component. We will focus on the cases when $s \in (d, 2d)$ in (1.3), in which percolation can be guaranteed, for instance, by requiring that the minimal probability of a nearest-neighbor connection, $p$, is sufficiently close to 1. (Indeed, in $d \geq 2$, it suffices that $p$ exceeds the percolation threshold for bond percolation on $\mathbb{Z}^d$, while in $d = 1$, this follows by the classic result of [22].) Moreover, by an extension of Burton–Keane's uniqueness argument due to [15], the infinite component is unique almost surely. We will use $\mathscr{C}_\infty$ to denote the set of sites in the infinite component of $\mathscr{G}$.

Our main result is as follows:

THEOREM 1.1. *Suppose that* (1.3) *holds with an* $s \in (d, 2d)$ *and assume that,* $\mathbb{P}$-*almost surely, the random graph* $\mathscr{G}$ *contains a unique infinite component* $\mathscr{C}_\infty$. *Then for all* $\varepsilon > 0$,

$$(1.6) \qquad \lim_{|x| \to \infty} \mathbb{P}\left(\Delta - \varepsilon \leq \frac{\log D(0, x)}{\log \log |x|} \leq \Delta + \varepsilon \,\Big|\, 0, x \in \mathscr{C}_\infty\right) = 1,$$

*where* $\Delta = \Delta(s, d)$ *is as in* (1.2).

Formula (1.6) is a precise form of the asymptotic expression (1.1). The fact that $\Delta^{-1}$ is the *binary* logarithm of $2d/s$ is a consequence of the fact that the longest bonds in the shortest path(s) between two distant sites of $\mathscr{C}_\infty$ exhibit a natural binary hierarchical structure; see Section 2 for more explanation. Note that $s \mapsto \Delta(s, d)$ is increasing throughout $(d, 2d)$ and, in particular, $\Delta(s, d) > 1$ for all $s \in (d, 2d)$ with $\lim_{s \downarrow d} \Delta(s, d) = 1$ and $\lim_{s \uparrow 2d} \Delta(s, d) = \infty$.

REMARK 1.1. The requirement of translation invariance is presumably not crucial for (the essence of) the above result. Indeed, most of our proofs should carry through under the weaker assumption of approximate homogeneity on large spatial scales. Notwithstanding, some of our arguments in Section 3 are based on previous results that require translation invariance and so we stick with the present setting for the rest of this paper.

1.3. *Discussion.* As already alluded to, several different asymptotic behaviors are possible in the above problem depending on the value of the exponent $s$. We proceed by reviewing the known (and conjectured) results. Throughout, we will focus on the specific distribution

$$(1.7) \qquad p_{xy} = 1 - \exp\{-\beta(1 + |x - y|)^{-s}\},$$

where $\beta \in [0, \infty)$. (Some of the results also required that all nearest-neighbor connection be a priori present.) We will concentrate on the asymptotic of two quantities: The *typical graph distance* $D(x, y)$—the focus of this paper—and the *diameter* $D_N$ of the graph obtained by "decorating" a box of $N \times \cdots \times N$ sites in $\mathbb{Z}^d$ by the bonds in $\mathscr{G}$ with both endpoints therein. There are five distinct regimes marked by the position of $s$ relative to the numbers $d$ and $2d$.

The cases of $s < d$ fall into the category of problems that can be analyzed using the concept of stochastic dimension introduced in [7]. The result is the almost-sure equality

$$(1.8) \qquad \sup_{x,y \in \mathbb{Z}^d} D(x, y) = \left\lceil \frac{d}{d - s} \right\rceil;$$

see Example 6.1 in [7]. A similar asymptotic statement holds for the $N \to \infty$ limit of $D_N$; see Theorem 4.1 of [6].

For $s = d$, Coppersmith, Gamarnik and Sviridenko [14] study the asymptotic of $D_N$. The resulting scaling is expressed by the formula

$$(1.9) \qquad D_N = \Theta(1)\frac{\log N}{\log \log N}, \qquad N \to \infty,$$

where $\Theta(1)$ is a quantity bounded away from 0 and $\infty$. Since the typical distance is always less than the diameter, this shows that $D(x, y)$ will grow at most logarithmically with $|x - y|$. However, at present the appropriate lower bound on $D(x, y)$ is missing.

In the cases $d < s < 2d$, Benjamini and Berger [6] and Coppersmith, Gamarnik and Sviridenko [14] proved polylogarithmic upper bounds on $D_N$ [and hence on $D(x, y)$ for $|x - y| \approx N$]. However, the best lower bound these references gave was proportional to $\log N$. The present paper provides a sharp leading-order asymptotic for $D(x, y)$ which constrains $D_N$ to grow at least as fast as $(\log N)^{\Delta + o(1)}$. Unfortunately, neither the bounds from [6] and [14] nor those derived for $D(x, y)$ in the present paper are sharp enough to make any definitive asymptotic statements about $D_N$. We hope to return to this question in a future publication.

The critical cases $s = 2d$ are at present not very well understood. Here Benjamini and Berger [6] conjectured that

$$(1.10) \qquad D_N = N^{\theta(\beta) + o(1)}, \qquad N \to \infty,$$

with $\theta(\beta) \in (0,1)$, and we expect a similar asymptotic to be valid for the typical distance $D(x,y)$. A general upper bound on the above $\theta(\beta)$ was derived in [14]. The corresponding—but not sharp—lower bounds were derived in [6] and [14] under the restriction to the "nonpercolative" regime $d = 1$ and $\beta < 1$. (This restriction appears because the proof relies heavily on the notion of a cut-point; see [22].) Surprisingly, similarly open are the cases $s > 2d$ where we expect that $D_N$ scales linearly with $N$. The latter seems to have been proved only in $d = 1$ [6], or for the case of the supercritical nearest-neighbor percolation in $d \geq 2$ [5]; see also [4].

An important technical resource for this paper has been the recent work of Berger [8] on long-range percolation with exponents $d < s < 2d$. Employing a variant of the renormalization scheme of [22], Berger proved among other things the *absence of critical percolation* and, whenever there *is* percolation, the existence of a cluster of at least $N^{d-o(1)}$ sites in any box of volume $N^d$ (see Theorem 3.1). An extension of this result (see Theorem 3.2) establishing tight bounds on the probability that the largest connected component in a finite box contains a *positive fraction* of all sites is essential for the proof of the upper bound in (1.6).

1.4. *Relation to "small-world" phenomena.* As already mentioned, long-range percolation has been recently used in the study of "small-world" phenomena. The catchy term "small worlds" originates in the old but still-fun-to-read article by Milgram [21], who observed through an ingenious experiment that two typical Americans are just six acquaintances (or six "handshakes") away from each other. With the rise of the overall world connectivity in recent years due to the massive expansion of air traffic, electronic communications and particularly the internet, and so on, the "small-world" phenomena experienced a fair amount of new interest. Novel examples emerged in physics and biology, particularly after the publication of [26]. Several mathematical models were devised and studied using both rigorous and nonrigorous techniques. A brief overview of the situation from the perspective of the theory of random graphs (and additional references) can be found in Section 10.6 of [9].

While we will not attempt to summarize the content of the publication boom following the appearance of [26], let us mention that a major deficiency of many models introduced so far seems to be—at least to the author of the present paper—the unclear status of their *relevance* to the actual (physical, biological or sociological) systems of interest. In particular, a large fraction of studied models seem to unjustly ignore the underlying spatial structure present in the practical problem of interest. (The reason for that is most likely the reduced complexity—as in statistical mechanics, models without underlying geometry, the so-called *mean-field* models, are often exactly solvable.) With this problem in mind, Benjamini and Berger [6] proposed a new class of "small-world" models based on long-range percolation on Euclidean graphs. More precisely, as an underlying graph they consider an a priori connected ring of $N$ sites to which long edges are added with probability as described in (1.7).

One of the questions discussed by Benjamini and Berger was how the diameter of the resulting random graph depends on $N$ for various ranges of values of $s$. As detailed in Section 1.3, this behavior depends rather sensitively on the value of the exponent $s$. In particular, "phase transitions" occur at $s = d$, which is the borderline of the region with finite diameters, and $s = 2d$, which separates the regions of linear and sublinear scaling. Each of the resulting behaviors may be useful in different contexts. For instance, if we believe Milgram's assertion that six is the typical graph distance between two average Americans regardless of the population size, the exponent $s$ should be within the regime described by (1.8).

**2. Main ideas of the proof.** The proof of Theorem 1.1 consists of two parts where we separately prove the upper and lower bounds in (1.6). Both parts will be based on the concept of certain *hierarchies* of sites whose definition is given below. In this definition—and elsewhere in this paper—the symbol $\sigma$ denotes a hierarchical index, $\sigma \in \{0,1\}^k$, which can be viewed as a parametrization of the leaves of a *binary* tree of depth $k$. Thus, for instance, $\sigma = 01101$ means that, starting from the root, we "go" left, right, right, left and right to reach the leaf represented by $\sigma$. Adding digits behind $\sigma$ denotes index concatenation.

DEFINITION 2.1. Given an integer $n \geq 1$ and distinct sites $x, y \in \mathbb{Z}^d$, we say that the collection

$$(2.1) \qquad \mathcal{H}_n(x,y) = \{(z_\sigma) : \sigma \in \{0,1\}^k, k = 1, 2, \ldots, n;\ z_\sigma \in \mathbb{Z}^d\}$$

is a *hierarchy of depth $n$ connecting $x$ and $y$* if:

1. $z_0 = x$ and $z_1 = y$.
2. $z_{\sigma 00} = z_{\sigma 0}$ and $z_{\sigma 11} = z_{\sigma 1}$ for all $k = 0, 1, \ldots, n-2$ and all $\sigma \in \{0,1\}^k$.
3. For all $k = 0, 1, \ldots, n-2$ and all $\sigma \in \{0,1\}^k$ such that $z_{\sigma 01} \neq z_{\sigma 10}$, the bond between $z_{\sigma 01}$ and $z_{\sigma 10}$ is occupied, that is, $(z_{\sigma 01}, z_{\sigma 10}) \in \mathcal{G}$.
4. Each bond $(z_{\sigma 01}, z_{\sigma 10})$ as specified in part 3 appears only once in $\mathcal{H}_k(x,y)$.

In the following, the pairs of sites $(z_{\sigma 00}, z_{\sigma 01})$ and $(z_{\sigma 10}, z_{\sigma 11})$ will be referred to as "gaps."

REMARK 2.1. By assumption 2, a hierarchy of depth $n$ is uniquely specified by its $n$th level. Note that we do not require the *sites* of the hierarchy to be distinct and, if two points of the form $z_{\sigma 10}$ and $z_{\sigma 11}$ coincide, we do not insist on having a bond between them. The phrase "connecting $x$ and $y$" in the definition of $\mathcal{H}_n(x, y)$ is not to imply that $\mathcal{H}_n(x, y)$ is an occupied path from $x$ to $y$. Instead, $\mathcal{H}_n(x, y)$ should be thought of as a framework of large-scale bonds which can be turned into a path by connecting the "gaps" in an appropriate way; see Figure 1.

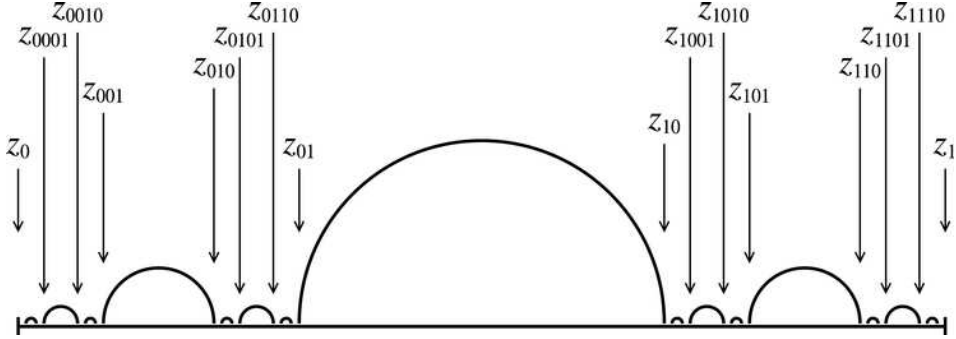

Fig. 1. *A schematic picture of a hierarchy of depth* 5 *connecting* $x = z_0$ *and* $y = z_1$. *The straight line represents a path between* $x$ *and* $y$ *and the arches stand for the bonds between pairs of sites* $(z_{\sigma 01}, z_{\sigma 10})$. *The arrows indicate the sites on levels* 1–4 *of the hierarchy; the fifth level has been omitted for brevity. Note that, by part* 2 *of Definition* 2.1, *we have* $z_{100} = z_{10}$, *and so on.*

Our strategy in both the upper and lower bound will be to identify a hierarchy of sufficient depth from within a given path. In an idealized situation, this hierarchy between sites at Euclidean distance $N$ would be such that the primary bond $(z_{01}, z_{10})$ has length (approximately) $N^{s/(2d)}$, the secondary bonds have length $N^{(s/(2d))^2}$, and so on. The principal difficulty is to "make" the hierarchy deep enough so that it already contains "most of" the bonds in the underlying path. In particular, we will have to guarantee that the "gaps"—which may still be rather spread out in the Euclidean distance—can be spanned without substantially affecting the overall length.

2.1. *Upper bound.* To outline the proof of the upper bound on the graph distance, it is convenient to start by analyzing the cases when all pairs of nearest neighbors on $\mathbb{Z}^d$ are a priori connected. In these situations one can (essentially) construct a path connecting two distant sites which uses about the optimal number of distinct occupied bonds.

Let $\gamma \in (s/(2d), 1)$. The construction is based on the following observation: If $x$ and $y$ are two sites at distance $|x - y| = N \gg 1$, and if $B_0$, respectively, $B_1$, are boxes of side $N^\gamma$ centered at $x$, respectively, $y$, then $B_0$ and $B_1$ are with overwhelming probability connected by a bond in $\mathcal{G}$. Indeed, there are $N^{d\gamma}$ ways to place each endpoint of such a bond while its Euclidean length has to be essentially equal to $N$. Hence, the probability that $B_0$ and $B_1$ are *not* directly connected by an occupied bond is

$$(2.2) \quad \mathbb{P}(B_0 \not\leftrightarrow B_1) = \exp\left\{ -\sum_{z \in B_0} \sum_{z' \in B_1} q(z - z') \right\} = \exp\{-N^{2d\gamma - s + o(1)}\}.$$

Since $2d\gamma > s$, the right-hand side tends rapidly to zero as $N \to \infty$.

Once the bond between $B_0$ and $B_1$ has been selected, a similar argument shows that the boxes $B_{00}$ and $B_{01}$ of side $N^{\gamma^2}$, centered at $x$ and the "nearer" endpoint of the primary bond, respectively, will typically be connected by an occupied bond. Continuing the process for a while, the endpoints of the family of bonds thus identified give rise to a hierarchy of sites in the above sense: First we let $z_0 = x$ and $z_1 = y$, then we let $z_{01}$ and $z_{10}$ be the endpoints of the primary bond connecting $B_0$ and $B_1$. Next, the endpoints of the secondary bonds connecting the boxes $B_{00}$ and $B_{01}$, respectively, $B_{10}$ and $B_{11}$, will be denoted by $z_{001}$ and $z_{010}$, respectively, $z_{101}$ and $z_{110}$. The higher levels will be denoted similarly. Note that, in order to have each level of the hierarchy completely defined, we need to use part 2 of Definition 2.1 to identify $z_{00}$ with $z_0$, and so on.

Of course, the most pertinent question is now for how long we can continue developing such a hierarchy. Proceeding as in our previous estimates, the probability that *not all* pairs of boxes $B_{\sigma 0}$ and $B_{\sigma 1}$ with $\sigma \in \{0, 1\}^k$ will be connected by a bond in $\mathscr{G}$ is bounded by

$$(2.3) \quad 2^k \exp\{-N^{\gamma^k(2d\gamma - s + o(1))}\},$$

where $2^k$ counts the number of bonds we are trying to control at this step and the factor $N^{\gamma^k(2d\gamma - s)}$ in the exponent originates from the fact that we are connecting boxes of side $N^{\gamma^{k+1}}$ which are at Euclidean distance $N^{\gamma^k}$ from each other. This estimate shows that, as long as $N^{\gamma^k(2d\gamma - s)} \gg \log \log N$, the probability that the identification procedure fails is negligible. However, this allows us to reach the level when the pairs of sites constituting the "gaps" are no farther than $N^{\gamma^k} = (\log N)^{o(1)}$ from each other. This happens for $k \approx K$, where

$$(2.4) \quad K = \frac{\log \log N}{\log(1/\gamma)}.$$

Now, a hierarchy of depth $K$ consists of roughly $2^K$ bonds and $2^K$ "gaps." Using nearest-neighbor paths to span each "gap," the total number of bonds needed to connect all "gaps" will thus be at most $2^K (\log N)^{o(1)}$. Hence, the graph distance between $x$ and $y$ cannot exceed $2^K (\log N)^{o(1)}$. Plugging the value of $K$ and passing to the limit $\gamma \downarrow s/(2d)$, the latter is no more than $(\log N)^{\Delta + o(1)}$.

Performing the above argument without the luxury of an a priori connected graph involves quite some extra work. Indeed, we need to ensure that the sites identified in the process are connected to $x$ and $y$ (and, therefore, to $\mathscr{C}_\infty$) and that the bonds lie in a "backbone"—rather than a "dead-end"—of the connection between $x$ and $y$. Our solution to this nonlocal optimization problem is to construct the hierarchy so that each site $z_\sigma$ for $\sigma \in \{0,1\}^k$ is connected to a positive fraction of all sites in the $N^{\gamma^k}$ neighborhood of $z_\sigma$. Since the distance between the endpoints of the "gaps" in such a hierarchy is at most of the order $N^{\gamma^{k-1}}$, the connected components of these endpoints are still with a large probability connected by a bond from $\mathscr{G}$. Now, if $k \approx K$, we have $N^{d\gamma^k} = (\log N)^{o(1)}$ and we need no more than $(\log N)^{o(1)}$ steps to connect the endpoints of each "gap." This allows us to proceed as before.

To ensure the connectivity property, we will introduce the concept of a *dense* site which is a site $x$ that is connected to at least a (prescribed) fraction of all sites in a sufficiently large box centered at $x$. Then we need to establish two additional facts: First, any site $x \in \mathscr{C}_\infty$ is with overwhelming probability dense. Second, any sufficiently large box contains a positive fraction of dense sites. These statements—which come as Corollaries 3.3 and 3.4—will allow us to look for hierarchies containing only dense sites, for which the above argument easily carries through. The proof of the two corollaries in turn requires showing that the largest connected component in any box contains a (uniformly) positive fraction of all sites. To maintain generality, this statement—which comes as Theorem 3.2—has to be proved under very modest assumptions; essentially, we only assume the asymptotic (1.3) and the fact that there is percolation.

2.2. *Lower bound.* The argument for the upper bound shows that there exists a path that connects $x$ to $y$ in about $(\log|x - y|)^\Delta$ steps. The goal of the lower bound is to show that, among the multitude of paths possibly connecting $x$ and $y$, no path will be substantially shorter.

In an idealized situation, our reasoning would go as follows: We set $N = |x - y|$ and pick a path from $\Pi(x, y)$ that connects $x$ with $y$ in less than $(\log N)^{O(1)}$ steps. [Here $O(1)$ represents a fixed number whose value is irrelevant in the following.] Next we will attempt to identify a hierarchy from $\pi$. The primary bond $(z_{01}, z_{10})$ is chosen simply as (one of) the longest bonds in $\pi$. Since $|x - y| = N$ but $|\pi| = (\log N)^{O(1)}$, this bond must be longer than $N/(\log N)^{O(1)}$. But in order for this bond to exist with a reasonable probability, a similar argument as used in the upper bound shows that the distances $N_0 = |x - z_{01}|$ and $N_1 = |z_{10} - y|$ must be such that

$$(2.5) \qquad\qquad N_0^d N_1^d \gtrsim N^{s + o(1)}.$$

Supposing (without any good reason) that $N_0$ is comparable with $N_1$, the removal of $(z_{01}, z_{10})$ from $\pi$ would leave us with two paths that connect sites at distance $N^{s/(2d)+o(1)}$ in a polylogarithmic number of steps. The argument could then be iterated which would eventually allow us to categorize the whole path into a hierarchical structure, with one bond of length $N$, two bonds of length $N^{s/(2d)}$, four bonds of length $N^{(s/(2d))^2}$, and so on.

It is easy to check that the hierarchy thus identified would involve roughly $2^K$ bonds, where $K$ is as in (2.4) with $\gamma = s/(2d)$, and $|\pi|$ would thus have to be at least $(\log N)^{\Delta+o(1)}$. Of course, the main problem with the above argument is that the assumption $N_0 \approx N_1$ is not justified and presumably fails in a large number of places. Extreme ways to violate the condition $N_0 \approx N_1$ are not so hard to dismiss. For instance, in the case of a "gap" collapse, for example, when $N_0 = N^{o(1)}$, the bound (2.5) forces that $N_1 \gtrsim N^{s/d+o(1)} \gg N$, implying that $(z_{01}, z_{10})$ was not the longest bond after all. But, since we are dealing with an exponentially growing number of bonds, even "soft" violations of this condition could make the whole argument crumble to pieces. As we will describe below, the solution is to work with (2.5)—and its generalizations—the way it stands without trying to extract any information about the particular $N_0$ and $N_1$.

Here is what we do. We pick a $\gamma$ satisfying $\gamma < s/(2d)$ and show that *every* path connecting $x$ and $y$ in less than $(\log N)^{O(1)}$ steps, where $N = |x - y|$, contains a hierarchy $\mathcal{H}_n(x, y)$ of depth $n \lesssim K$—with $K$ as in (2.4)—such that the following holds with overwhelming probability: The length of the "gaps" is comparable with the length of the bonds that "fill" them; that is, for all $k = 1, \ldots, n-1$, we have

$$(2.6) \qquad |z_{\sigma 01} - z_{\sigma 10}| = |z_{\sigma 0} - z_{\sigma 1}|^{1-o(1)}, \qquad \sigma \in \{0, 1\}^{k-1}.$$

Moreover, the *average* size of the "gaps" on the $k$th level is about $N^{\gamma^k}$; that is, the quantities $N_\sigma = |z_{\sigma 0} - z_{\sigma 1}|$ satisfy

$$(2.7) \qquad \prod_{\sigma \in \{0, 1\}^k} N_\sigma \geq N^{(2\gamma)^k}$$

for all $k = 1, \ldots, n-1$. Obviously, for $k = 1$ this is a more precise form of (2.5). Part (2.6) is a consequence of the fact that, in order to connect two sites at distance $N^{\gamma^k}$ in less than $(\log N)^{O(1)}$ steps, at least one bond in the path must be longer than $N^{\gamma^k}/(\log N)^{O(1)}$. This equals $N^{\gamma^k(1-o(1))}$ as long as $k \lesssim K$. As to the proof of (2.7), let $\mathcal{E}_n$ be the event that the inequality in (2.7) holds for $k = 1, 2 \ldots, n-1$. We will sketch the derivation of an upper bound on $\mathbb{P}(\mathcal{E}_{k+1}^c \cap \mathcal{E}_k)$ which can then be iterated into a bound on $\mathbb{P}(\mathcal{E}_n^c)$.

Fix a collection of numbers $(N_\sigma)$ representing the distances between various "gaps" in the hierarchy, and let us estimate the probability that a hierarchy with these $(N_\sigma)$ occurs. In light of (2.6), the primary bond will cost $N^{-s+o(1)}$ of probability, but there are of the order $(N_0 N_1)^{d-1}$ ways to choose the endpoints. (Remember that $N_0$ and $N_1$ are fixed.) Similarly, the two secondary bonds cost $(N_0 N_1)^{-s+o(1)}$ of probability and their endpoints contribute of the order $(N_{00} N_{01} N_{10} N_{11})^{d-1}$ of entropy. Applying this to the collections $(N_\sigma)$ compatible with $\mathcal{E}_{n+1}^c \cap \mathcal{E}_n$, we get

$$
\begin{aligned}
(2.8) \qquad & \mathbb{P}(\mathcal{E}_{n+1}^c \cap \mathcal{E}_n) \\
& \leq \sum_{(N_\sigma)} \frac{(N_0 N_1)^{d-1}}{N^{s-o(1)}} \frac{(N_{00} N_{01} N_{10} N_{11})^{d-1}}{(N_0 N_1)^{s-o(1)}} \cdots \frac{\prod_{\sigma \in \{0,1\}^n} N_\sigma^{d-1}}{\prod_{\sigma \in \{0,1\}^{n-1}} N_\sigma^{s-o(1)}},
\end{aligned}
$$

where the sum goes over all $(N_\sigma)$ for which $\mathcal{E}_{n+1}^c \cap \mathcal{E}_n$ holds.

To evaluate the right-hand side of (2.8), we need to observe that the numerator $(N_0 N_1)^{d-1}$ can be combined with the denominator of the next quotient into a term which by $s > d$ is summable on both $N_0$ and $N_1$; using (2.7), the resulting sum over $N_0$ and $N_1$ is bounded by $N^{-(s-d)(2\gamma)+o(1)}$. The other numerators will be handled analogously; the upshot is that, for all $k \leq n-1$, the sum over all $N_\sigma$ with $\sigma \in \{0,1\}^k$ is bounded by $N^{-(s-d)(2\gamma)^k+o(1)}$. The last numerator has no denominator to be matched with, but here we can use that, since (2.7) for $k = n$ *fails* on $\mathcal{E}_{n+1}^c \cap \mathcal{E}_n$, the product of $N_\sigma$'s is now bounded from above! Consequently, the relevant sum does not exceed $N^{d(2\gamma)^n+o(1)}$. Putting all these estimates together, and applying the inequality

$$
(2.9) \qquad s - d(2\gamma)^n + (s-d) \sum_{k=1}^{n-1} (2\gamma)^k \geq (s - 2d\gamma)(2\gamma)^{n-1},
$$

the right-hand side of (2.8) is bounded by $N^{-(s-2d\gamma)(2\gamma)^{n-1}}$. [The inequality in (2.9) can be derived either by direct computation or by a repeated application of the inequality $s + (s-d)(2\gamma) \geq s(2\gamma)$ to the first and the third term on the left-hand side.] Summing the obtained bound over $n$, the probability $\mathbb{P}(\mathcal{E}_n)$ is shown to be essentially one as long as $n \lesssim K$.

Once we have established that (2.6) and (2.7) hold, we will use a similar estimate to find a lower bound on $D(x, y)$. Here we simply have to prove that, even though the hierarchy is already rather large, the lower bound (2.7) requires that at least as many bonds be used to connect all of the "gaps." To avoid some unpleasant combinatorial estimates, we will continue under the simplifying assumption that all of the $2^{n-1}$ "gaps" of the hierarchy are nontrivial.

Let $\overline{\mathcal{F}}_n$ be the event that every hierarchy of depth $n$ satisfying (2.6) and (2.7) requires more than $2^{n-1}$ extra steps to connect all of its "gaps." In light of our bound on $\mathbb{P}(\mathcal{E}_n^c)$, it suffices to estimate the probability of $\overline{\mathcal{F}}_n^c \cap \mathcal{E}_n$. Since all "gaps" are nontrivial, the only way $\overline{\mathcal{F}}_n^c$ can occur is that each "gap" is spanned by a single bond. Now the bond spanning the "gap" $(z_{\sigma 0}, z_{\sigma 1})$ costs $N_\sigma^{-s+o(1)}$ amount of probability and so $\mathbb{P}(\overline{\mathcal{F}}_n^c \cap \mathcal{E}_n)$ can be bounded by

$$
\begin{aligned}
(2.10) \quad & \mathbb{P}(\overline{\mathcal{F}}_n^c \cap \mathcal{E}_n) \\
& \leq \sum_{(N_\sigma)} \frac{(N_0 N_1)^{d-1}}{N^{s-o(1)}} \cdots \frac{\prod_{\sigma \in \{0,1\}^{n-1}} N_\sigma^{d-1}}{\prod_{\sigma \in \{0,1\}^{n-2}} N_\sigma^{s-o(1)}} \prod_{\sigma \in \{0,1\}^{n-1}} \frac{1}{N_\sigma^{s-o(1)}},
\end{aligned}
$$

where the $(N_\sigma)$'s now obey (2.7) for *all* $k = 0, 1, \ldots, n-1$. The last product on the right-hand side makes the entire sum convergent and (for $N \gg 1$) small. Thus, with overwhelming probability, $\overline{\mathcal{F}}_n$ occurs for all $n \lesssim K$, which means that the shortest path(s) between $x$ and $y$ must contain at least $2^{n-1} = 2^{K(1-o(1))}$ distinct bonds. For $K$ as in (2.4), we have $2^K = (\log N)^{\Delta'}$, where $1/\Delta' = \log_2(1/\gamma)$. From here the lower bound in (1.6) follows by letting $\gamma \uparrow s/(2d)$.

## 3. Percolation in finite boxes.

3.1. *Size of the largest connected component.* In this section we will prove an estimate showing that the largest connected component in large but finite boxes contains a positive fraction of all sites whenever there is percolation. This estimate will be essential for the proof of the upper bound in (1.6). Throughout this section, the original meaning of the quantity from (1.3) will be substituted by a weaker form (3.1) below. We will return to the original definition in Section 4.

We begin by quoting a result from [8]. Let us say that the collection of probabilities $(p_{xy})_{x,y \in \mathbb{Z}^d}$ is *percolating*, if the associated i.i.d. measure has an infinite cluster almost surely. Let $L \geq 1$ be an integer and let $\Lambda_L$ be a box in $\mathbb{Z}^d$ of side $L$ containing $L^d$ sites. Consider the percolation problem restricted to the sites of $\Lambda_L$ (and, of course, only the bonds with both endpoints in $\Lambda_L$) and let $|\mathscr{C}_L|$ denote the size of the largest connected component in $\Lambda_L$. In [8], Berger proved that once $(p_{xy})_{x,y \in \mathbb{Z}^d}$ are percolating, $\Lambda_L$ contains a large cluster. The precise formulation is as follows:

THEOREM 3.1 ([8], Lemma 2.3). *Let $d \geq 1$ and suppose that the collection of probabilities $(p_{xy})_{x,y \in \mathbb{Z}^d}$, where $p_{xy} = 1 - e^{-q(x-y)}$, is percolating. Suppose that, for some $s \in (d, 2d)$,*

$$
(3.1) \qquad \liminf_{|x| \to \infty} |x|^s (1 - e^{-q(x)}) > 0.
$$

*Then for each $\varepsilon > 0$ and each $\zeta \in (0, \infty)$, there exists an $L$ such that*

$$
(3.2) \qquad \mathbb{P}(|\mathscr{C}_L| < \zeta L^{s/2}) \leq \varepsilon.
$$

We note that once (3.1) holds for some $s$, then it holds also for any $s' > s$. Therefore, Theorem 3.1 actually guarantees that the largest connected component in $\Lambda_L$ will contain at least $L^{d-o(1)}$ sites. (Note that the statement forbids us to take $s = 2d$, and an inspection of Lemma 2.3 in [8] reveals that this is nontrivially rooted in the proof.) However, for our purposes we need to work with the event that $|\mathscr{C}_L|$ is *proportional* to $L^d$ and, in addition, we also need a more explicit estimate on the probability of such an event. Our extension of Berger's result comes in the following theorem:

THEOREM 3.2.   *Let $d \geq 1$ and consider the probabilities $(p_{xy})_{x,y \in \mathbb{Z}^d}$ such that (3.1) holds for some $s \in (d, 2d)$. Suppose that $(p_{xy})_{x,y \in \mathbb{Z}^d}$ are percolating. For each $s' \in (s, 2d)$ there exist numbers $\rho > 0$ and $L_0 < \infty$ such that for each $L \geq L_0$,*

$$(3.3) \qquad \mathbb{P}(|\mathscr{C}_L| < \rho|\Lambda_L|) \leq e^{-\rho L^{2d-s'}}.$$

*In particular, once $L$ is sufficiently large, the largest connected component in $\Lambda_L$ typically contains a positive fraction of all sites in $\Lambda_L$.*

Theorem 3.2 alone would allow us to establish the existence of a hierarchy between two sites, but it would not ensure that the "gaps" are properly connected (which is what we need to turn the hierarchy into a path). Fortunately, the structure of the proof of Theorem 3.2 allows us to make this conclusion anyway. To state the relevant mathematical claims, for each $x \in \mathbb{Z}^d$ and any odd integer $L \geq 1$, let $\Lambda_L(x)$ be the box of side $L$ centered at $x$ and let $\mathscr{C}_L(x)$ be the set of sites in $\Lambda_L(x)$ that are connected to $x$ by an occupied path in $\Lambda_L(x)$. Then we have:

COROLLARY 3.3.   *Under the conditions of Theorem 3.2, there exists a constant $\rho > 0$ such that*

$$(3.4) \qquad \lim_{L \to \infty} \mathbb{P}(|\mathscr{C}_L(x)| < \rho|\Lambda_L(x)|, x \in \mathscr{C}_\infty) = 0$$

*holds for each $x \in \mathbb{Z}^d$.*

COROLLARY 3.4.   *Given $\ell < L$, let $\mathscr{D}_L^{(\rho,\ell)}$ be the set of sites $x \in \Lambda_L$ such that $|\mathscr{C}_\ell(x)| \geq \rho|\Lambda_\ell(x)|$. Under the conditions of Theorem 3.2, for each $s' \in (s, 2d)$ there are constants $\ell_0 < \infty$ and $\rho > 0$ such that*

$$(3.5) \qquad \mathbb{P}(|\mathscr{D}_L^{(\rho,\ell)}| < \rho|\Lambda_L|) \leq e^{-\rho L^{2d-s'}}$$

*holds for any $\ell$ with $\ell_0 \leq \ell \leq L/\ell_0$.*

Theorem 3.2 and the two corollaries are what this section contributes to the proof of the main result of this paper. An impatient (or otherwise uninterested) reader should feel free to skip the rest of this section on a first reading and pass directly to Sections 4.1 and 4.2. For those who stay put, we proceed by discussing the main ideas of the proof and a breakdown of its steps into the corresponding technical lemmas. The actual proofs appear in Section 3.5.

3.2. *Outline of the proof.*  Our strategy of the proof of Theorem 3.2 is as follows. First we combine a one-step renormalization with Theorem 3.1 to convert the problem into a similar question for *site-bond* percolation. An important feature of this reformulation is that the occupation probability of sites and bonds can be made as close to 1 as we wish.

Given an odd integer $K \geq 1$, let $\mathscr{C}_K^{(x)}$ denote the largest connected component in $\Lambda_K(x)$; in the case of a tie we will choose the component containing the site that is minimal in the standard lexicographic order on $\mathbb{Z}^d$. For any two distinct $x, y \in \mathbb{Z}^d$, we will say that $\Lambda_K(Kx)$ and $\Lambda_K(Ky)$ are *directly connected* if there is an occupied bond connecting a site from $\mathscr{C}_K^{(Kx)}$ to a site from $\mathscr{C}_K^{(Ky)}$. We will use $\{\Lambda_K(Kx) \nleftrightarrow \Lambda_K(Ky)\}$ to denote the event that $\Lambda_K(Kx)$ and $\Lambda_K(Ky)$ are not directly connected. Then we have:

LEMMA 3.5.  *Under the assumptions of Theorem 3.1, for any $s \in (d, 2d)$ the following is true: For each $\beta < \infty$ and $r < 1$ there exist a number $\delta > 0$ and an odd integer $K < \infty$, such that*

$$(3.6) \qquad \mathbb{P}(|\mathscr{C}_K^{(Kx)}| < \delta \, |\Lambda_K(Kx)|) \leq 1 - r$$

*and*

$$(3.7) \qquad \mathbb{P}(\Lambda_K(Kx) \nleftrightarrow \Lambda_K(Ky)) \leq e^{-\beta |x-y|^{-s}}$$

*hold for all distinct $x, y \in \mathbb{Z}^d$.*

Regarding boxes of side $K$ as new sites and the pairs of maximal-connected components connected by a bond from $\mathscr{G}$ as new bonds, Lemma 3.5 allows us to set up a renormalization scheme of [22]. Clearly, by (3.6) and (3.7), a site is occupied with probability at least $r$ and two occupied sites $x$ and $y$ are connected with probability at least $1 - e^{-\beta |x-y|^{-s}}$. (For the sites that are not occupied such a connection will not be relevant, so we will often assume that the latter holds for all sites.) This puts us into a position where we can apply the following "idealized" version of the desired claim:

LEMMA 3.6.  *Let $d \geq 1$ and consider the site-bond percolation model on $\mathbb{Z}^d$ with sites being occupied with probability $r \in [0, 1]$ and the bond between sites $x$ and $y$ being occupied with probability*

$$(3.8) \qquad p_{xy} = 1 - \exp\{-\beta |x-y|^{-s}\},$$

*where $s \in (d, 2d)$ and $\beta \geq 0$. Let $|\mathscr{C}_N|$ denote the size of the largest connected component of occupied sites and occupied bonds in $\Lambda_N$. For each $s' \in (s, 2d)$ there exist numbers $N_0 < \infty$, $\vartheta > 0$ and $\beta_0 < \infty$ such that*

$$(3.9) \qquad \mathbb{P}_{\beta, r}(|\mathscr{C}_N| < \vartheta|\Lambda_N|) \leq e^{-\vartheta \beta N^{2d-s'}}$$

*holds true for all $N \geq N_0$ whenever $\beta \geq \beta_0$ and $r \geq 1 - e^{-\vartheta \beta}$.*

REMARK 3.1.  The fact that the exponent in (3.9) is proportional to $\beta$ will not be needed for the proofs in this paper. The addition of $\beta$ represents only a minor modification of the proof and the stronger result will (hopefully) facilitate later reference.

Once Lemma 3.6 is inserted into the game, the proof of Theorem 3.2 will be easily concluded. To prove Lemma 3.6, we will invoke a combination of coarse-graining, stochastic domination and a corresponding estimate for the *complete graph*. (Let us recall that a complete graph of $n$ vertices is a graph containing a bond for each unordered pair of distinct numbers from $\{1, 2, \ldots, n\}$.) The relevant complete-graph statement is extracted into the following lemma:

LEMMA 3.7.  *Consider a complete graph of $n$ vertices. Let each site be occupied with probability $r$ and each bond be occupied with probability $p$. Let $\mathbb{P}_n^{p,r}$ be the resulting i.i.d. measure and let $|\mathscr{C}_n|$ denote the number of sites in the largest connected component of occupied sites and occupied bonds. For each $q, q' \in [0, 1]$ with $q' < q$, there exists a number $\psi(q', q) > 0$ such that for each $r' \in [0, r)$, each $p' \in [0, p)$ and all $n \geq 1$,*

$$(3.10) \qquad \mathbb{P}_n^{p,r}(|\mathscr{C}_n| \leq p'r'n) \leq e^{-n\psi(r',r)} + e^{-(1/2)(n^2(r')^2-n)\psi(p',p)}.$$

*Moreover, for each $\alpha \in [0, 1)$, there exists a constant $C = C(\alpha) < \infty$ such that*

$$(3.11) \qquad \psi(q', q) \geq (1-\alpha)(1-q')\left[\log \frac{1}{1-q} - C\right]$$

*holds true for all $q, q' \in [0, 1]$ satisfying the bound $(1-q') \geq (1-q)^\alpha$.*

REMARK 3.2.  While Lemma 3.6 can presumably be proved without reference to the complete graph, in our case the passage through Lemma 3.7 has the advantage of easily obtained quantitative estimates. As mentioned before, the present proof of Theorem 3.2 invokes a renormalization scheme for which Theorem 3.1—whose proof is based on a similar renormalization scheme—serves as a starting point. It would be of some conceptual interest to see whether a more direct proof of Theorem 3.2 based on a single renormalization is possible.

REMARK 3.3.    In $d \geq 2$, the decay rate in (3.9)—and, consequently, in (3.3)—is not always optimal. The reason is that, for $\beta \gg 1$ and $1 - r \ll 1$, the site-bond percolation problem dominates the nearest-neighbor percolation on $\mathbb{Z}^d$ for which it is expected (and essentially proved, see [3, 13]) that the probability in (3.9) should decay exponentially with $N^{d-1}$. [For $s \in (d, 2d)$, this is sometimes better and sometimes worse than $N^{2d-s}$.] This alternative decay rate is not reflected in our proofs because, to apply equally well in all dimensions $d \geq 1$, they consistently rely only on long-range connections.

Having outlined the main strategy and stated the principal lemmas, we can plunge into the proofs. First we will prove technical Lemmas 3.5 and 3.7. Then we will establish the site-bond percolation Lemma 3.6. Once all preparatory statements have been dispensed with, we will assemble the ingredients into the proofs of Theorem 3.2 and Corollaries 3.3 and 3.4.

3.3. *Preparations.*    Here we will prove Lemmas 3.5 and 3.7. First we will attend to the one-step renormalization from Lemma 3.5 whose purpose is to wash out the short-range irregularities of the $p_{xy}$'s and to ensure that the constants $\beta$ and $r$ in Lemma 3.6 can be chosen as large as required.

PROOF OF LEMMA 3.5.    The principal task before us is to choose $K$ so large that both bounds follow from Theorem 3.1 and the assumption (3.1). Let $d \geq 1$ and let $\beta < \infty$ and $r < 1$. Let $\varepsilon = 1 - r$ and pick a $d' \in (s/2, d)$ such that $d' - (d - 1) > d - d'$. By Theorem 3.1 and the paragraph that follows, for each $\zeta > 0$, there exists a $K = K(\varepsilon, d', \zeta)$ such that $|\mathscr{C}_K| \geq \zeta K^{d'}$ occurs with probability exceeding $1 - \varepsilon$. Since Theorem 3.1 allows us to make $K$ arbitrarily large (indeed, the constraint $|\mathscr{C}_K| \leq K^d$ forces $K$ to exceed a positive power of $\zeta$), we can also assume that, for some constant $\alpha > 0$,

$$(3.12) \qquad p_{xy} \geq 1 - e^{-\alpha|x-y|^{-s}} \qquad \text{once } |x - y| \geq K^{d-d'}.$$

Here we rely on (3.1).

Let $b > 0$ be a constant such that $bK$ exceeds the diameter of $\Lambda_K$ in the $|\cdot|$-norm for all $K$. We will show that (3.6) and (3.7) hold once $K$ is large enough and, in particular, so large that

$$(3.13) \qquad \frac{1}{4} \frac{\alpha \zeta^2}{(2b)^s} K^{2d'-s} \geq \beta.$$

(Note that $2d' > s$ so the left-hand side increases with $K$.) Consider a partitioning of $\mathbb{Z}^d$ into disjoint boxes of side $K$; that is, let us write $\mathbb{Z}^d = \bigcup_{x \in \mathbb{Z}^d} \Lambda_K(Kx)$. We will call a box $\Lambda_K(Kx)$ occupied if the bond configuration restricted to $\Lambda_K(Kx)$ contains a connected component larger than $\zeta K^{d'}$. By the choice of $K$, each $\Lambda_K(Kx)$ is occupied independently with probability exceeding $1 - \varepsilon = r$. This proves (3.6) with $\delta = \zeta K^{d'-d}$.

To prove also (3.7), we need to ensure that sufficient portions of the components in the two $K$-blocks are so far from each other that (3.12) can safely be applied. To this end, we note that, since $|\mathscr{C}_K^{(Kx)}| \geq \zeta K^{d'}$, at least half of the sites in $\mathscr{C}_K^{(Kx)}$ will be farther from $\mathbb{Z}^d \setminus \Lambda_K(Kx)$ than $\eta_K = a\zeta K^{d'-(d-1)}$, where $a$ is a constant depending only on the norm $|\cdot|$ and the dimension. Moreover, an easy argument shows that if $x, y \in \mathbb{Z}^d$ are distinct and $z \in \mathscr{C}_K^{(Kx)}$ and $z' \in \mathscr{C}_K^{(Ky)}$ are such that $\mathrm{dist}(z, \mathbb{Z}^d \setminus \Lambda_K(Kx)) \geq \eta_K$ and $\mathrm{dist}(z', \mathbb{Z}^d \setminus \Lambda_K(Ky)) \geq \eta_K$, then

$$(3.14) \qquad \eta_K |x - y| \leq |z - z'| \leq 2bK|x - y|.$$

Now our choice of $d'$ guarantees that, if $K$ is sufficiently large, we will have $\eta_K \geq K^{d-d'}$ and the inequality on the left-hand side shows that (3.12) for $p_{zz'}$ is in power. The bound on the right-hand side of (3.14) then allows us to write

$$(3.15) \quad \mathbb{P}(\Lambda_K(Kx) \leftrightarrow \Lambda_K(Ky)) \leq \exp\left\{-\alpha \frac{(\zeta K^{d'}/2)^2}{(|x - y|2bK)^s}\right\} \leq e^{-\beta|x-y|^{-s}},$$

where we used (3.13) to derive the last inequality. This completes the proof. $\square$

Next we will focus on the proof of Lemma 3.7, which concerns the complete graph. This lemma will be used to drive the induction argument in the next section.

PROOF OF LEMMA 3.7.   The proof starts by estimating the total number of occupied sites. Once that number is known to be sufficiently large, the desired bound is a result of conditioning on occupied sites combined with straightforward estimates concerning occupied bonds.

Fix $r' \in (0, r)$ and $p' \in [0, p)$, and let $\rho = p'r'$. Note that we can assume that $r'n > 1$, because otherwise the right-hand side of (3.10) exceeds 1. Let $A_n$ denote the (random) number of occupied vertices of the complete graph. Since $A_n$ can be represented as a sum of independent random variables with mean $r$, the event $\{A_n < r'n\}$ forces $A_n$ to deviate substantially from its mean and the standard Chernoff bound implies

$$(3.16) \qquad \mathbb{P}_n^{p,r}(A_n < r'n) \leq e^{-\psi(r',r)n},$$

where

$$(3.17) \qquad \psi(q', q) = \sup_{\lambda \geq 0}[-\log(1 - q + qe^{-\lambda}) - \lambda q'].$$

As is easy to show, $\psi(q', q) > 0$ for all $q' < q$.

The bound (3.16) is responsible for the first term on the right-hand side of (3.10). It remains to show that the conditional probability given $A_n \geq r'n$ is bounded by the second term. Thus suppose that $A_n \geq r'n$ and let $V_n$ denote the total number of unordered pairs of occupied sites *not* connected by an occupied bond. The principal observation is that if $|\mathscr{C}_n| \leq \rho n$, then $V_n$ has to be rather large. More precisely, we claim that $|\mathscr{C}_n| \leq \rho n$ implies $V_n \geq \frac{1}{2}(A_n^2 - \rho n A_n)$. Indeed, if we label the connected components of occupied sites and bonds by index $i$ and use $k_i$ to denote the number of sites in the $i$th occupied component, then the number of vacant bonds certainly exceeds

$$(3.18) \qquad \sum_{i<j} k_i k_j = \tfrac{1}{2}\left(\sum_i k_i\right)^2 - \tfrac{1}{2}\sum_i k_i^2.$$

On $\{|\mathscr{C}_n| \leq \rho n\}$ we have $k_i \leq \rho n$ for each $i$ and since also $\sum_i k_i = A_n$, the second sum can be bounded by $A_n \rho n$. The desired inequality $V_n \geq \frac{1}{2}(A_n^2 - \rho n A_n)$ follows.

In light of our previous reasoning, we are down to estimating the probability

$$(3.19) \qquad \mathbb{P}_n^{p,r}(V_n \geq \tfrac{1}{2}(A_n^2 - \rho n A_n)|A_n \geq r'n).$$

The estimate will be performed by conditioning on the set of occupied sites. Once the set of occupied sites has been fixed, $V_n$ can be represented as a sum of $N = \frac{1}{2}(A_n^2 - A_n)$ independent random variables, each of which has mean $1 - p$. Now, assuming $A_n \geq r'n$ and recalling that $\rho = p'r' \leq r'$ and $r'n > 1$, we can estimate

$$(3.20) \quad \tfrac{1}{2}(A_n^2 - \rho n A_n) = N\left(1 - \frac{\rho n - 1}{A_n - 1}\right) \geq N\left(1 - \frac{\rho n - 1}{r'n - 1}\right) \geq N(1 - p').$$

(In the cases when $\rho n < 1$, we just skip the intermediate inequality.) Now, since $p' < p$, the event $\{V_n \geq \frac{1}{2}(A_n^2 - \rho n A_n)\}$ constitutes a large deviation for the random variable $V_n$. Invoking again the Chernoff bound and a little algebra, we thus get

$$(3.21) \quad \mathbb{P}_n^{p,r}(V_n \geq \tfrac{1}{2}(A_n^2 - \rho n A_n)|A_n \geq r'n, \tfrac{1}{2}(A_n^2 - A_n) = N) \leq e^{-N\psi(p',p)}.$$

From here (3.10) follows by noting that on $\{A_n \geq r'n\}$ we have $N \geq \frac{1}{2}(n^2(r')^2 - n)$.

To verify the bound (3.11), we need to find the minimizing $\lambda$ in (3.17) and use it to find an explicit expression for $\psi(q',q)$. A computation gives

$$(3.22) \qquad \psi(q',q) = (1 - q')\log\left(\frac{1 - q'}{1 - q}\right) - q'\log\left(\frac{q}{q'}\right).$$

The second term—including the minus sign—can be split into two parts: the term $q'\log q'$, which is bounded below by $-1/e$, and the term $-q'\log q$, which is always positive. Moreover, if $q' \to 1$, then $(q'\log q')/(1 - q') \to 1$. From here we infer that the second term is bounded below by $(1 - \alpha)(1 - q')$-times a (negative) $\alpha$-dependent constant. Using the bound $1 - q' \geq (1 - q)^\alpha$ in the first term, (3.11) is proved.  □

3.4. *Site-bond percolation lemma.* Now we are ready to start proving Lemma 3.6. The essential part of the proof relies on induction along a series of scales defined as follows. Fix an $s' \in (s, 2d)$, let $\ell_0$ be a positive integer and consider an increasing sequence $(\ell_n)$ of integers such that $\ell_n \geq \ell_0$ for all $n \geq 1$. Let $N_0$ be another positive integer and let

$$(3.23) \qquad N_n = N_0 \prod_{k=1}^{n} \ell_k.$$

Suppose that $\ell_n$ tend to infinity so fast that

$$(3.24) \qquad \sum_{n \geq 1} \ell_n^{d-s} < \infty \quad \text{and} \quad \sum_{n \geq 1} N_n^{s-s'} < \infty$$

but moderate enough that also

$$(3.25) \qquad c_0 = \inf_{n \geq 1} \frac{1}{\ell_{n+1}^{s'}} \prod_{k=1}^{n} \ell_k^{2d-s'} > 0.$$

Next, let us define the sequences $(r_n)$ and $(p_n)$ by putting $r_n = 1 - 6\ell_n^{d-s}$ and $p_n = 1 - N_n^{s-s'}$, let us fix a constant $\rho_0 \in (0, 1)$ and let

$$(3.26) \qquad \rho_n = \rho_0 \prod_{k=1}^{n} (r_k p_k).$$

Clearly, by (3.24) we have $\rho_n \to \rho_\infty > 0$ as $n \to \infty$.

REMARK 3.4. An extreme example of a sequence $(\ell_n)$ satisfying these constraints is

$$(3.27) \qquad \ell_n = e^{(1+a)^{n-1}}, \qquad n \geq 1,$$

where $a = (2d - s')/s'$. Here $a$ has been tuned in such a way that the term in the infimum in (3.25) is independent of $n$. Since $a > 0$, the bounds (3.24) immediately follow.

The proof of Lemma 3.6 is based on the fact that a bound (3.9) for $N = N_n$ can be used to prove the same bound for (essentially) any $N$ between $N_n$ and $N_{n+1}$. The proof works as soon as $c_0$ and $\rho_\infty$ are bounded away from zero and $\ell_0$ and $N_0$ are sufficiently large; the precise form of $(\ell_n)$ is not important. [In particular, we do not make *explicit* use of (3.24).] The induction step is isolated into the following claim:

LEMMA 3.8. *Suppose the assumptions of Lemma 3.6 hold and let $s' \in (s, 2d)$, $c_0' > 0$ and $\tau_0 > 0$. Then there exist two numbers $\ell_0' < \infty$ and $N_0' < \infty$ and a constant $c_1 \in (0, \infty)$ such that for any $N_0 \geq N_0'$ and any sequence $(\ell_n)$* *with $\ell_0 \geq \ell_0'$ and $c_0 \geq c_0'$, the following holds: If $\tau \in [\tau_0, c_1 \rho_\infty^2 \beta]$ and if $k$ is a nonnegative integer such that*

$$
\mathbb{P}_{\beta,r}(|\mathscr{C}_{N_k}| \leq \rho_k N_k^d) \leq e^{-\tau N_k^{2d-s'}}
\tag{3.28}
$$

*is true, then*

$$
\mathbb{P}_{\beta,r}(|\mathscr{C}_{\ell N_k}| \leq \rho_{k+1}(\ell N_k)^d) \leq e^{-\tau(\ell N_k)^{2d-s'}}
\tag{3.29}
$$

*holds for all $\ell \in \{\ell_0, \ldots, \ell_{k+1}\}$.*

PROOF.   The proof is based on the "complete-graph" Lemma 3.7. Our preliminary task is to set up all the constants so that the bounds emerging from this lemma are later easily converted to that on the right-hand side of (3.29). Let $s' \in (s, 2d)$, $c_0' > 0$ and $\tau_0 > 0$ be fixed. First we will address the choice of the constants $\ell_0'$, $N_0'$ and $c_1$. We will assume that $\ell_0'$ is so large that

$$
(1 - 6\ell^{d-s})^2 \ell^{2d} - \ell^d \geq \tfrac{1}{2}\ell^{2d}
\tag{3.30}
$$

holds for all $\ell \geq \ell_0'$. Then we choose $N_0'$ so large that, for all $N_0 \geq N_0'$ and all $\tau \geq \tau_0$,

$$
6\ell^{d-s} > e^{-(1/2)\tau(\ell N_0)^{2d-s'}}
\tag{3.31}
$$

holds for all $\ell \geq \ell_0'$ and that

$$
\tau N_0^{2d-s'} \geq \max\{3C, \log 2\},
\tag{3.32}
$$

where $C$ is the constant from Lemma 3.7 for $\alpha = 1/2$. Moreover, we let $b$ be the constant such that $Nb$ bounds the diameter of $\Lambda_N$ in the metric $|\cdot|$ for all $N \geq 1$ and choose $c_1 = b^{-s}/16$. Then we also require that $N_0'$ be so large that for any $N_0 \geq N_0'$,

$$
\beta b^{-s} \rho_n^2 N_n^{2d-s} \ell^{-s} - C \geq 8\tau N_n^{2d-s} \ell^{-s}
\tag{3.33}
$$

and

$$
(\ell N_n)^{s-s'} > e^{-(1/2)\beta b^{-s} \rho_n^2 N_n^{2d-s} \ell^{-s}}
\tag{3.34}
$$

hold for all $\tau \leq c_1 \rho_\infty^2 \beta$, all $n \geq 0$ and all $\ell \in \{\ell_0, \ldots, \ell_{n+1}\}$. The first bound is verified for sufficiently large $N_0$ by noting the inequalities $\beta b^{-s} \rho_n^2 \geq 2\tau_0$ and

$$
N_n^{2d-s} \ell^{-s} \geq N_n^{2d-s} \ell_{n+1}^{-s} \geq N_0^{2d-s}\left(\ell_{n+1}^{-s'} \prod_{k \leq n} \ell_k^{2d-s'}\right)^{s/s'} \geq N_0^{2d-s} c_0^{s/s'}.
\tag{3.35}
$$

As to (3.34), we note that $N_n^{2d-s} \ell^{-s}/(\ell N_n)^{s'-s} \geq N_n^{2d-s'} \ell_{n+1}^{-s'}$, which again can be made arbitrarily large by boosting up $N_0$. The factors $\tfrac{1}{2}$ in the exponents of (3.31) and (3.34) have been inserted for later convenience.

Now we are ready to prove (3.29) from (3.28). Let $\ell_0 \geq \ell'_0$ and $N_0 \geq N'_0$. Suppose that (3.28) holds for some $\tau \in [\tau_0, c_1 \rho_\infty^2 \beta]$ and some $k \geq 0$. Pick an $\ell \in \{\ell_0, \ldots, \ell_{k+1}\}$. Viewing $\Lambda_{\ell N_k}$ as the disjoint union of $\ell^d$ translates $\Lambda^{(i)}$ of $\Lambda_{N_k}$,

$$(3.36) \qquad \Lambda_{\ell N_k} = \bigcup_{i=1}^{\ell^d} \Lambda^{(i)},$$

let us call $\Lambda^{(i)}$ occupied if it contains a connected component of size at least $\rho_k |N_k|^d$. Choosing one such connected component in each occupied $\Lambda^{(i)}$ (if necessary, employing lexicographic order on $\mathbb{Z}^d$ to resolve potential ties), we will call $\Lambda^{(i)}$ and $\Lambda^{(j)}$ with $i \neq j$ connected if their respective connected components are directly connected by an occupied bond. Let

$$(3.37) \qquad r = 1 - e^{-\tau N_k^{2d-s'}} \quad \text{and} \quad p = 1 - e^{-\beta b^{-s} \rho_k^2 N_k^{2d-s} \ell^{-s}}$$

and note that (3.28) implies that $r$ is a lower bound on the probability that $\Lambda^{(i)}$ is occupied. Moreover, a simple calculation shows that $p$ is a uniform lower bound on the probability that two distinct $\Lambda^{(i)}$ and $\Lambda^{(j)}$ are connected.

Let us introduce the quantities

$$(3.38) \qquad r' = 1 - 6\ell^{d-s} \quad \text{and} \quad p' = 1 - (\ell N_k)^{s-s'}$$

and let $\mathcal{A}_{k,\ell}$ be the event that the largest connected component of the occupied $\Lambda^{(i)}$'s comprises more than $p'r'\ell^d$ individual boxes. The assumption (3.31) shows that $r' < r$, while (3.34) guarantees that $p' < p$. Invoking the fact that $\mathcal{A}_{k,\ell}$ is an increasing event, the probability of $\mathcal{A}_{k,\ell}^c$ is bounded by the probability that, for site-bond percolation on a complete graph with $\ell^d$ vertices and parameters $r$ and $p$, the largest connected component involves at most $p'r'\ell^d$ vertices:

$$(3.39) \qquad \mathbb{P}_{\beta,r}(\mathcal{A}_{k,\ell}^c) \leq \mathbb{P}_{\ell^d}^{p,r}(|\mathscr{C}_{\ell^d}| \leq p'r'\ell^d).$$

Since the factors $\frac{1}{2}$ in the exponents in (3.31) and (3.34) ensure that $1 - r' \geq (1-r)^{1/2}$ and $1 - p' \geq (1-p)^{1/2}$, the right-hand side can be estimated using Lemma 3.6 with $\alpha = \frac{1}{2}$ and $C = C(\frac{1}{2})$. To evaluate the first term on the right-hand side of (3.10), we estimate

$$(3.40) \qquad \ell^d \psi(r', r) \geq \tfrac{1}{2} 6\ell^{2d-s}(\tau N_k^{2d-s'} - C) \geq 2\tau(\ell N_k)^{2d-s'}.$$

Here we used that, by our choice of $\alpha$, we have $\ell^d(1-\alpha)(1-r') = \frac{1}{2}6\ell^{2d-s}$ and then invoked (3.32) to show that $\tau N_k^{2d-s'} \geq \tau N_0^{2d-s'} \geq 3C$. Similarly, we get

$$(3.41) \qquad \begin{aligned} \tfrac{1}{2}(\ell^{2d}(r')^2 - \ell^d)\psi(p', p) &\geq \tfrac{1}{4}\ell^{2d}(\ell N_k)^{s-s'}(\beta b^{-s} \rho_k^2 N_k^{2d-s} \ell^{-s} - C) \\ &\geq 2\tau(\ell N_k)^{2d-s'} \end{aligned}$$

for the exponent in the second term in (3.10). Here we first used (3.30) to reduce the complicated $\ell$-dependence on the extreme left and then we inserted (3.11) and the definitions of $p'$ and $p$ to produce the intermediate inequality. Finally, we invoked (3.33).

By putting the bounds (3.10), (3.11) and (3.40), (3.41) together and recalling that $\log 2 \leq \tau N_0^{2d-s'}$, the probability $\mathbb{P}_{\beta,r}(\mathcal{A}_{k,\ell}^c)$ does not exceed the term on the right-hand side of (3.29). But on $\mathcal{A}_{k,\ell}$, the box $\Lambda_{\ell N_k}$ contains a connected component comprising (strictly) more than $p'r'\ell^d$ disjoint connected components, each of which involves at least $\rho_k N_k^d$ sites. Using that $r' \geq r_{k+1}$ and $p' \geq p_{k+1}$, we have

$$(3.42) \qquad |\mathcal{C}_{\ell N_k}| > p'r'\ell^d \rho_k |N_k|^d \geq \rho_{k+1}(\ell N_k)^d \qquad \text{on } \mathcal{A}_{k,\ell},$$

and thus $\{|\mathcal{C}_{\ell N_k}| \leq \rho_{k+1}|\Lambda_{\ell N_k}|\} \subset \mathcal{A}_{k,\ell}^c$. From here (3.28) follows. $\quad\square$

Lemma 3.8 encapsulates the induction step. However, we will also need an estimate that allows us to start the induction. This is provided in the following lemma.

LEMMA 3.9. *Under the conditions of Lemma 3.6, for each $c_2 \in (0,\infty)$, there exist numbers $N_0 < \infty$ and $\vartheta_0 \in (0,1)$ and, for each $\vartheta < \vartheta_0$, there exists a number $\beta_0 < \infty$ such that*

$$(3.43) \qquad \mathbb{P}_{\beta,r}\big(|\mathcal{C}_N| < \tfrac{1}{4}|\Lambda_N|\big) \leq e^{-c_2\vartheta\beta N^{2d-s}}$$

*holds once $N \geq N_0$, $\beta \geq \beta_0$ and $r \geq 1 - e^{-\vartheta\beta}$.*

PROOF. We will again apply the "complete-graph" Lemma 3.7. Let $p' = r' = \frac{1}{2}$, let $b$ be a constant such that $bN$ exceeds the diameter of $\Lambda_N$ for any $N$ and pick $\vartheta_0 > 0$ such that $128c_2\vartheta_0 < b^{-s}$. Fix a number $\vartheta \in (0,\vartheta_0)$. Then the left-hand side of (3.43) is bounded by the left-hand side of (3.10) with $n = N^d$, $p = 1 - \exp\{-\beta(bN)^{-s}\}$ and $r = 1 - e^{-\vartheta\beta}$. We will estimate the right-hand side of (3.10) under the conditions when $N_0$ is so large that

$$(3.44) \qquad N^{s-d} \geq 16c_2 \quad \text{and} \quad N^{2d}(r')^2 - N \geq N^{2d}/8$$

are true for all $N \geq N_0$ and, given such an $N$, the constant $\beta_0$ is so large that for all $\beta \geq \beta_0$ we have $1 - r' \geq (1-r)^{1/2}$ and $1 - p' \geq (1-p)^{1/2}$, and

$$(3.45) \qquad \vartheta\beta \geq 2C \quad \text{and} \quad \beta(bN)^{-s} - C \geq 128\vartheta\beta N^{-s}.$$

Here, as before, $C$ is the constant from Lemma 3.7 for $\alpha = \frac{1}{2}$.

In conjunction with these bounds, (3.11) with $\alpha = \frac{1}{2}$ shows that

$$(3.46) \qquad \psi(r',r) \geq \tfrac{1}{4}(\vartheta\beta - C) \geq \tfrac{1}{8}\vartheta\beta$$

and

$$(3.47) \qquad \psi(p', p) \geq \tfrac{1}{4}(\beta(bN)^{-s} - C) \geq 32c_2 \vartheta \beta N^{-s}.$$

Using the bounds in (3.44), we find that both exponents on the right-hand side of (3.10) exceed $2c_2 \vartheta \beta N^{2d-s}$. This implies $\mathbb{P}_{\beta,r}(|\mathscr{C}_N| \leq \vartheta |\Lambda_N|) \leq 2 \exp\{-2c_2 \vartheta \times \beta N^{2d-s}\}$. Increasing $\beta_0$ if necessary, the latter is no more than $\exp\{-c_2 \vartheta \beta N^{2d-s}\}$.

$\square$

Equipped with the induction machinery from Lemmas 3.8 and 3.9, the proof of the main site-bond percolation lemma is now easily concluded.

PROOF OF LEMMA 3.6.  First we will adjust the parameters so that Lemmas 3.8 and 3.9 can directly be applied. Let $s' \in (s, 2d)$ and let $c_0' > 0$ and $\tau_0 > 0$ be fixed. Let $\ell_0'$, $N_0'$ and $c_1$ be the constants from Lemma 3.8 and pick a sequence $(\ell_0)$ such that $\ell_0 \geq \ell_0'$ and $c_0 \geq c_0'$ are satisfied. Pick a number $N_0 \geq N_0'$ so large that Lemma 3.9 holds for $c_2 = (2\ell_0)^{2d-s'}$ and $N \geq N_0$, and let $\vartheta_0$ be the corresponding constant from this lemma. Let $\rho_0 \in (0, \tfrac{1}{4}]$ and define $\rho_n$ and $\rho_\infty$ as in (3.26). Let $\vartheta > 0$ be such that $\vartheta \leq \vartheta_0$, $(2\ell_0)^d \vartheta \leq \rho_\infty$ and $c_2 \vartheta \leq c_1 \rho_\infty^2$. Choose $\beta_0$ so large that Lemma 3.9 holds for all $\beta \geq \beta_0$ and such that $c_2 \vartheta \beta_0 > \tau_0$. Note that $\tau = c_2 \vartheta \beta$ necessarily satisfies $\tau \in [\tau_0, c_1 \rho_\infty^2 \beta]$ as long as $\beta \geq \beta_0$, which is needed in Lemma 3.8.

Now we are ready to run the induction argument: Since $\rho_0 \leq \tfrac{1}{4}$ and $\vartheta \leq \vartheta_0$, Lemma 3.9 ensures that (3.28) holds for $\tau = c_2 \vartheta \beta$ and $k = 0$. Applying the induction step from Lemma 3.8, we recursively show that (3.29) is true for all $k \geq 0$ and all $\ell \in \{\ell_0, \ldots, \ell_{k+1}\}$. Let $N$ be a general integer and let $k$ be a nonnegative integer such that $N_{k+1} > N \geq N_k$. Let $\ell$ be the maximal number in $\{1, \ell_0, \ldots, \ell_{k+1}\}$ such that $\ell N_k \leq N$. A simple calculation now shows that $|\Lambda_N| \leq \max(2\ell_0)^d |\Lambda_{\ell N_k}|$ and, if we position $\Lambda_N$ and $\Lambda_{\ell N_k}$ so that $\Lambda_{\ell N_k} \subset \Lambda_N$,

$$(3.48) \qquad \{|\mathscr{C}_N| \leq \vartheta |\Lambda_N|\} \subset \{|\mathscr{C}_{\ell N_k}| \leq (2\ell_0)^d \vartheta |\Lambda_{\ell N_k}|\}.$$

By our previous conclusions and the fact that $(2\ell_0)^d \vartheta \leq \rho_{k+1}$, the probability of the event on the right-hand side is bounded by $\exp\{-\tau(\ell N_k)^{2d-s'}\}$. From here, (3.9) for a general $N$ follows by noting that, by our choice of $c_2$, we have $\tau(\ell N_k)^{2d-s'} \geq \vartheta \beta N^{2d-s'}$.  $\square$

3.5. *Proofs of Theorem 3.2 and Corollaries 3.3 and 3.4*.  Now we are finally ready to prove Theorem 3.2. After some preliminary arguments, the proof follows a line of reasoning similar to the one just used to prove Lemma 3.6.

PROOF OF THEOREM 3.2. Let $s' \in (s, 2d)$ and let $\vartheta$, $\beta_0$ and $N_0$ be as in Lemma 3.6. Pick numbers $\beta \geq \beta_0$ and $r \geq 1 - e^{-\vartheta\beta}$ and let $K$ and $\delta$ be the corresponding constants from Lemma 3.5. First we will prove the claim for $L$ of the form $L = KN$, where $N$ is a positive integer. To that end, let us view $\Lambda_L$ as the disjoint union of $N^d$ translates $\Lambda_K(Kx)$ of $\Lambda_K$, where $x \in \mathbb{Z}^d$. We will call $\Lambda_K(Kx)$ occupied if $|\mathscr{C}_K^{(Kx)}| \geq \delta|\Lambda_K|$ and, similarly, we will call two distinct $\Lambda_K(Kx)$ and $\Lambda_K(Ky)$ connected if the connected components $\mathscr{C}_K^{(Kx)}$ and $\mathscr{C}_K^{(Ky)}$—chosen with the help of lexicographic order in case of a tie—are directly connected by a bond from $\mathscr{G}$.

By Lemma 3.5 and our choice of $\delta$ and $K$, the box $\Lambda_K(x)$ is occupied with probability exceeding $r$, while $\Lambda_K(Kx)$ and $\Lambda_K(Ky)$ are connected with probability exceeding $p_{xy}$ in (3.8). Let $\mathcal{A}_{N,K}$ be the event that the box $\Lambda_L$ contains a connected component $\mathscr{C}_N$ of boxes $\Lambda_K(Kx)$ such that at least $\vartheta|\Lambda_N|$ of the connected components in these boxes get joined in $\mathscr{C}_N$. By Lemma 3.6, we know that

$$\mathbb{P}_{\beta,r}(\mathcal{A}_{N,K}^{c}) \leq e^{-\vartheta\beta N^{2d-s'}}. \tag{3.49}$$

On the other hand, on $\mathcal{A}_{N,K}$ we have

$$|\mathscr{C}_L| \geq (\vartheta|\Lambda_N|)(\delta|\Lambda_K|) = \vartheta\delta|\Lambda_L|, \tag{3.50}$$

and thus $\{|\mathscr{C}_L| \leq \rho|\Lambda_L|\} \subset \mathcal{A}_{N,K}^{c}$ once $\rho < \vartheta\delta$. If $\rho$ is also less than $\vartheta\beta K^{s'-2d}$, this finishes the proof for $L$ of the form $NK$. The general values of $L$ are handled by noting that if $NK \leq L < (N+1)K$, then $|\Lambda_L| \leq 2^d|\Lambda_{NK}|$ and, if $\Lambda_{NK} \subset \Lambda_L$, then also $|\mathscr{C}_{NK}| \leq |\mathscr{C}_L|$.   $\square$

PROOF OF COROLLARY 3.3. Fix $s' \in (s, 2d)$, let $N_0$, $\vartheta$ and $\beta_0$ be the constants from Lemma 3.6, and let $\beta \geq \beta_0$ and $r \geq 1 - e^{-\vartheta\beta}$. Let $\ell > 3N_0$ be an odd integer and let $K$ be the constant from Lemma 3.5 for our choice of $\beta$ and $r$. Clearly, it suffices to show that (3.4) holds for $L$ of the form $L_n = K\ell^n$ and $\rho$ proportional to the product of constants $\delta$ and $\vartheta$ from Lemmas 3.5 and 3.6. All of the volumes $\Lambda_{L_n}$ below are centered at $x$ so we omit that fact from the notation.

Our strategy is as follows: We pick an $\varepsilon > 0$ and show that, with probability at least $1 - \varepsilon$ and some integer $n$, the largest connected component $\mathscr{C}_{L_n}$ in $\Lambda_{L_n}$ is connected to at least $\rho|\Lambda_{L_{n'}}|$ sites in $\Lambda_{L_{n'}}$, for every $n' \geq n$. (Note that this guarantees that $\mathscr{C}_{L_n} \subset \mathscr{C}_\infty$.) Once this has been established, we observe that $|\mathscr{C}_{L_{n'}}(x)| < \rho|\Lambda_{L_{n'}}|$ implies that $x$ cannot be connected to $\mathscr{C}_{L_n}$ within $\Lambda_{L_{n'}}$. Assuming that $x \in \mathscr{C}_\infty$, the box $\Lambda_{L_n}$ then contains at least two distinct sites $x, y \in \mathscr{C}_\infty$ which are not connected within $\Lambda_{L_{n'}}$. By the uniqueness of the infinite cluster, the probability of the latter event can be made smaller than $\varepsilon$ by making $n'$ sufficiently large. But then the limit in (3.4) must be less than $2\varepsilon$ and, since $\varepsilon$ was arbitrary, it must equal to zero.

To make the proof complete, it remains to establish the first claim in the previous paragraph. Namely, we must show that the probability that $\mathscr{C}_{L_n}$ is *not* connected to at least $\rho|\Lambda_{L_{n'}}|$ sites in $\Lambda_{L_{n'}}$ for *some* $n' \geq n$ is less than $\varepsilon$, provided $n$ is sufficiently large. To that end, let $B_k$, with $k \geq 0$, be a sequence of boxes (generally not centered at $x$) such that $B_0 = \Lambda_{L_n}$ and that $B_k$ is the maximal box in $\Lambda_{L_{n+k}}$ that is centered on the $(1, 1, \ldots, 1)$ half-axis, disjoint from all the previous $B_k$'s and with side a multiple of $K$. Since $\ell > 3$, it is easy to see that $|B_k|$ grows proportionally to $|\Lambda_{L_{n+k}}|$ (in fact, $|B_k|/|\Lambda_{L_{n+k}}| \geq 3^{-d}$). Our goal is to show that the largest connected components in all $B_k$'s are with overwhelming probability connected.

Invoking Lemmas 3.5 and 3.6 and choosing $n$ sufficiently large, the probability that each box $B_k$—viewed as the disjoint union of translates of $\Lambda_K$—contains a component comprising at least $\vartheta|B_k|/K^d$ maximal connected components $\mathscr{C}_K^{(Kx)}$ of size at least $\delta K^d$ is bounded by $\varepsilon/2$. On the other hand, the probability that the corresponding components in $B_k$ and $B_{k+1}$ are *not* connected is bounded by

$$(3.51) \qquad \exp\left\{-\beta\frac{\vartheta^2}{K^{2d}}\frac{|B_k||B_{k+1}|}{(bL_{n+k+1}/K)^s}\right\} \leq \exp\{-\beta'\ell^{(2d-s)(n+k)}\},$$

where $b$ is a constant such that $bL_n/K$ bounds the distance $|x - y|$ for any translates $\Lambda_K(Kx)$ and $\Lambda_K(Ky)$ contained in $\Lambda_{L_n}$, and where $\beta'/\beta$ is a constant that depends only on $\ell$ and $K$. The right-hand side is summable on $k$ and the sum can be made smaller than $\varepsilon/2$ by increasing $n$. Thus with probability as least $1 - \varepsilon$, for each $k \geq 0$ the component $\mathscr{C}_{L_n}$ is connected to at least $\delta\vartheta|B_k|$ sites in $\Lambda_{L_{n+k}}$. Choosing $\rho < \delta\vartheta 3^{-d}$, the above claim follows. □

PROOF OF COROLLARY 3.4. Let $\ell > 1$ be an odd integer. Clearly, it suffices to prove the result for $\mathscr{D}_L^{(\rho,2\ell)}$ instead of $\mathscr{D}_L^{(\rho,\ell)}$ and $L$ a multiple of $3\ell$. Viewing $\Lambda_L$ as a disjoint union of boxes $\Lambda_{3\ell}(3\ell x)$ with $x \in \mathbb{Z}^d$, let $\mathscr{C}_\ell^{(3\ell x)}$ be a maximal connected component in $\Lambda_\ell(3\ell x)$. For $\rho' > 0$, let $\mathcal{A}_\ell(x)$ be the event that $|\mathscr{C}_\ell^{(3\ell x)}| \geq \rho'|\Lambda_\ell|$ is true. Let $N_{L,\ell}$ denote the number of $x \in \mathbb{Z}^d$ with $\Lambda_{3\ell}(3\ell x) \subset \Lambda_L$ such that $\mathcal{A}_\ell(x)$ occurs.

The events $\mathcal{A}_\ell(x)$ are independent and, if $\ell$ is large enough and $\rho' > 0$ is sufficiently small, Theorem 3.2 shows that $\mathcal{A}_\ell(x)$ occurs with probability at least $r = 1 - \exp\{-\rho'\ell^{2d-s'}\}$. This allows us to invoke the Chernoff bound once again with the result

$$(3.52) \qquad \mathbb{P}(N_{L,\ell} < r'L^d/(3\ell)^d) \leq e^{-\psi(r',r)L^d/(3\ell)^d},$$

where $\psi(r', r)$ is as in (3.17). Choosing $r' = \frac{1}{2}$, $\alpha = \frac{1}{2}$ and $C = C(\frac{1}{2})$ and taking $\ell$ so large that $\rho'\ell^{2d-s'} \geq 2C$, (3.11) gives us

$$(3.53) \qquad \psi(r', r) \geq \frac{1}{4}(\rho'\ell^{2d-s'} - C) \geq \frac{1}{8}\rho'\ell^{2d-s'}.$$

On the other hand, on $\{N_{L,\ell} \geq r'L^d/(3\ell)^d\}$ we have

$$(3.54) \qquad |\mathscr{D}_L^{(\rho,2\ell)}| \geq (\rho'\ell^d)(r'L^d/(3\ell)^d) = \rho'r'3^{-d}L^d,$$

and so $\{|\mathscr{D}_L^{(\rho,2\ell)}| < \rho L^d\} \subset \{N_{L,\ell} < r'L^d/(3\ell)^d\}$ once $\rho \leq \rho'r'3^{-d}$. Invoking (3.52) and (3.53) and choosing $\rho$ such that also $\rho \leq \frac{1}{8}3^{-d}\rho'$, the desired estimate follows. $\quad\square$

## 4. Proof of main result.

In this section we will provide the proof of Theorem 1.1. The arguments closely follow the outline presented in Sections 2.1 and 2.2. The reader may consider skimming through these sections once again before plunging into the ultimate details of the proofs.

### 4.1. Upper bound.

The principal goal of this section is to establish the upper bound in (1.6). By the continuity of $s \mapsto \Delta(s, d)$ it suffices to prove the upper bound for any number $s$ exceeding the limit (1.3), so we will instead assume that $s$ obeys (3.1). The desired claim is then formulated as follows:

PROPOSITION 4.1. Let $s \in (d, 2d)$ be such that (3.1) holds and let $\Delta = \Delta(s, d)$ be as in (1.2). For each $\Delta' > \Delta$ and each $\varepsilon > 0$, there exists an $N_0 < \infty$ such that

$$(4.1) \qquad \mathbb{P}(D(x,y) \geq (\log|x-y|)^{\Delta'}, x, y \in \mathscr{C}_\infty) \leq \varepsilon$$

holds for all $x, y \in \mathbb{Z}^d$ with $|x - y| \geq N_0$.

As discussed in Section 2.1, the proof is conceptually rather simple: For each pair of sites $x$ and $y$ we will construct a hierarchy of an appropriate depth $k$ connecting $x$ and $y$, such that pairs $(z_{\sigma01}, z_{\sigma10})$ with $\sigma \in \{0,1\}^{k-2}$ are connected by paths of length $(\log|x-y|)^{o(1)}$. The main difficulty stems from the requirement that the bonds constituting the hierarchy be connected in a prescribed (linear) order. This will be ensured by the condition that all sites constituting the hierarchy are surrounded by a sufficiently dense connected component.

Recall our notation that $\Lambda_L(x)$ is a box of side $L$ centered at $x$ and $\mathscr{C}_L(x)$ is the set of sites in $\Lambda_L(x)$ connected to $x$ by a path in $\Lambda_L(x)$. We will require that the sites $z_\sigma$ are dense points according to the following definition:

DEFINITION 4.1. Given a number $\rho \in (0,1)$ and an odd integer $\ell > 1$, we will call $x \in \mathbb{Z}^d$ a $(\rho, \ell)$-dense (or, simply, dense) point if $|\mathscr{C}_\ell(x)| \geq \rho|\Lambda_\ell(x)|$.

For any real $L > 0$ sufficiently large, let $L^+$ be the minimal odd integer larger than $L$ and let $L^-$ be the minimal odd integer larger than $L/2$. Let

$$(4.2) \qquad B_L(x) = \Lambda_{L^+}(x) \setminus \Lambda_{L^-}(x).$$

Given a number $\rho \in (0,1)$ and an odd integer $\ell > 1$, let $\mathscr{D}_L^{(\rho,\ell)}(x)$ denote the set of all $(\rho,\ell)$-dense points in $B_L(x)$. The input needed from Section 3 then comes directly from Corollaries 3.3 and 3.4. By Corollary 3.4 and the fact that $B_L(x)$ contains a box of side at least $L/3$, we know that

$$\mathbb{P}(|\mathscr{D}_L^{(\rho,\ell)}(x)| \leq \rho L^d) \leq e^{-\rho L^{2d-s'}} \tag{4.3}$$

once $\rho$ is sufficiently small, $s' \in (s,2d)$ and $\ell_0 \leq \ell \leq L/\ell_0$. Corollary 3.3 in turn shows that if $x, y \in \mathscr{C}_\infty$, then both $x$ and $y$ are dense points in the sense that for each $\varepsilon > 0$ there exists an $\ell_1 = \ell_1(\varepsilon) < \infty$ such that

$$\mathbb{P}(|\mathscr{C}_\ell(x)| \leq \rho \ell^d, x \in \mathscr{C}_\infty) \leq \varepsilon \tag{4.4}$$

is true whenever $\ell \geq \ell_1$. A similar statement holds for $y$.

Now we can define the principal events: Let $\gamma \in (s/(2d), 1)$ and let $x$ and $y$ be two sites in $\mathbb{Z}^d$. Let $N = |x - y|$ and define $N_n = N^{\gamma^n}$. For each $n \geq 1$, let $\mathcal{B}_n = \mathcal{B}_{n,\gamma}^{(\rho,\ell)}(x,y)$ be the event that there exists a hierarchy $\mathcal{H}_n(x,y)$ of depth $n$ connecting $x$ and $y$ subject to the following constraints: For all $k = 0, 1, \ldots, n-2$ and all $\sigma \in \{0,1\}^k$,

$$z_{\sigma 01} \in \mathscr{D}_{N_{k+1}}^{(\rho,\ell)}(z_{\sigma 0}) \quad \text{and} \quad z_{\sigma 10} \in \mathscr{D}_{N_{k+1}}^{(\rho,\ell)}(z_{\sigma 1}). \tag{4.5}$$

The event $\mathcal{B}_n$ ensures that all sites of the hierarchy—*except perhaps $x$ and $y$*—are $(\rho,\ell)$-dense points. To cover these exceptions, we also introduce the event $\mathcal{T} = \mathcal{T}^{(\rho,\ell)}(x,y)$ that both $x$ and $y$ are $(\rho,\ell)$-dense in the above sense. In the following, we will regard the number $\rho$ as fixed—such that (4.3) and (4.4) hold—but $\ell$ and $\gamma$ will be adjustable.

The requirements (4.5) become appreciated in the proof of the following bound:

LEMMA 4.2.   *For each $\varepsilon \in (0,1)$, each $\gamma \in (s/(2d), 1)$ and each $\Delta'$ satisfying*

$$\Delta' > \frac{\log 2}{\log(1/\gamma)}, \tag{4.6}$$

*there exists a constant $N' = N'(\varepsilon, \gamma, \Delta') < \infty$ such that the following is true for all $x, y \in \mathbb{Z}^d$ with $N = |x-y| \geq N'$: Let $n$ be the maximal positive integer such that*

$$n \log(1/\gamma) \leq \log \log N - \varepsilon \log \log \log N. \tag{4.7}$$

*If $\ell$ in the definition of the events $\mathcal{B}_n = \mathcal{B}_{n,\gamma}^{(\rho,\ell)}(x,y)$ and $\mathcal{T} = \mathcal{T}^{(\rho,\ell)}(x,y)$ is an odd integer between $N_n$ and $2N_n$, then*

$$\mathbb{P}(\{D(x,y) \geq (\log N)^{\Delta'}\} \cap \mathcal{B}_n \cap \mathcal{T}) \leq \varepsilon. \tag{4.8}$$

The reason why we choose $n$ as in (4.7) can be seen from the following bounds:

$$(4.9) \quad 2^n \leq (\log N)^{\log 2/\log(1/\gamma)} \quad \text{and} \quad e^{(1/\gamma)(\log \log N)^\varepsilon} \geq N^{\gamma^n} \geq e^{(\log \log N)^\varepsilon}.$$

These bounds will be important in the upcoming proof.

PROOF OF LEMMA 4.2.    The main reason why $N$ has to exceed a certain constant is because we need the scales corresponding to successive levels of the hierarchy to be clearly separated. To that end we observe that for all $k \leq n$ and $N_k = N^{\gamma^k}$, we have

$$(4.10) \qquad \qquad \log \frac{N_k}{N_{k+1}} \geq (1-\gamma)(\log \log N)^\varepsilon,$$

which tends to infinity as $N \to \infty$.

Introduce the abbreviation $\overline{\mathcal{B}}_n = \mathcal{T} \cap \mathcal{B}_n$. If $\overline{\mathcal{B}}_n$ occurs, then there exists at least one hierarchy $\mathcal{H}_n(x, y)$ of depth $n$ connecting $x$ and $y$ such that (4.5) is satisfied. Then (4.10) guarantees the existence of numbers $N' < \infty$ and $b \in (0, 1)$ such that the following is true for any such hierarchy: If $\sigma \in \{0, 1\}^k$ with $k = 0, 1, \ldots, n-2$, then

$$(4.11) \qquad \qquad |z_{\sigma 01} - z_{\sigma 10}| \geq b N_k,$$

while if $z \in B_{N_{k+2}}(z_{\sigma 0})$ and $z' \in B_{N_{k+2}}(z_{\sigma 01})$, then

$$(4.12) \qquad \qquad b^{-1} N_{k+1} \geq |z - z'| \geq b N_{k+1},$$

whenever $N \geq N'$. Similar statements hold for pairs $z \in B_{N_{k+2}}(z_{\sigma 1})$ and $z' \in B_{N_{k+2}}(z_{\sigma 10})$. Moreover, $N'$ can be chosen so large that also the bounds

$$(4.13) \qquad b^{-1} N_{n-1} < b N_{n-2} \quad \text{and} \quad b N_{n-1} > \text{diam} \, \Lambda_\ell$$

hold true for all $\ell$ between $\ell_0$ and $2N_n$ and all $n$ satisfying (4.7).

Let $\mathcal{A}_n$ be the event that, for any hierarchy that would make $\overline{\mathcal{B}}_n$ satisfied, at least one of the "gaps" of the hierarchy, say $(z_{\sigma 0}, z_{\sigma 01})$ where $\sigma \in \{0, 1\}^{n-2}$, fails to have the components $\mathscr{C}_\ell(z_{\sigma 0})$ and $\mathscr{C}_\ell(z_{\sigma 01})$ connected by an occupied bond. (Note that these components are quite large because both $z_{\sigma 0}$ and $z_{\sigma 01}$ are dense points.) We claim that

$$(4.14) \qquad \{D(x, y) \geq (\log N)^{\Delta'}\} \cap \overline{\mathcal{B}}_n \subset \mathcal{A}_n \cap \overline{\mathcal{B}}_n.$$

Indeed, if all "gaps" *do have* the corresponding components connected, then each $z_{\sigma 0}$ is connected to $z_{\sigma 01}$ by a path of no more than $1 + 2\ell^d$ bonds [note that $\ell^d$ bonds should be enough to get out of $\mathscr{C}_\ell(z_{\sigma 0})$], and similarly for the pairs $z_{\sigma 1}$ and $z_{\sigma 10}$. Noting that a hierarchy of depth $n$ involves only $2^{n-1}$ "gaps" and $2^{n-1} - 1$ bonds, we can use $\ell \leq 2N_n$ and (4.9) to write

$$(4.15) \qquad \begin{aligned} D(x, y) &\leq 2^{n-1}(1 + 2^{d+1} N_n^d) + 2^{n-2} \\ &\leq 2^{d+2}(\log N)^{\log 2/\log(1/\gamma)} e^{(d/\gamma)(\log \log N)^\varepsilon}. \end{aligned}$$

In light of (4.6), this is not compatible with $D(x, y) \geq (\log N)^{\Delta'}$ if $N$ is sufficiently large.

To finish the proof, we thus need to estimate the probability of $\mathcal{A}_n \cap \overline{\mathcal{B}}_n$. The above estimates show that the occurrence of $\mathcal{B}_n$ is determined by looking only at the bonds longer than $bN_{n-2}$ (to ensure the existence of a hierarchy) or shorter than $\operatorname{diam} \Lambda_\ell$ [to ensure that all sites in the hierarchy are $(\rho, \ell)$-dense]. Explicitly, let $\mathscr{F}$ denote the $\sigma$-algebra generated by the random variables $(\omega_{zz'})$ with $|z - z'| \geq bN_{n-2}$ or $|z - z'| \leq \operatorname{diam} \Lambda_\ell$. Then (4.11) and (4.13) show that $\overline{\mathcal{B}}_n \in \mathscr{F}$. This allows us to prove (4.8) by conditioning: Let $\beta$ be a number such that $p_{zz'} \geq 1 - \exp\{-\beta|z - z'|^{-s}\}$ for any pair $(z, z')$ of sites with $|z - z'| \geq N_n$. Then we have

$$(4.16) \qquad \mathbb{P}(\mathcal{A}_n | \mathscr{F}) \leq 2^{n-1} \exp\left\{-\beta \rho^2 \frac{N_n^{2d}}{(N_{n-1}/b)^s}\right\} \qquad \text{on } \overline{\mathcal{B}}_n.$$

Here we used that, on $\overline{\mathcal{B}}_n$, the components $\mathscr{C}_\ell(z_{\sigma 0})$ and $\mathscr{C}_\ell(z_{\sigma 01})$ are both larger than $\rho\ell^d \geq \rho N_n^d$, while (4.12) dictates that the longest bond that can connect them is not longer than $N_{n-1}/b$. The prefactor represents the number of "gaps" in the hierarchy, which is the number of places where $\mathcal{A}_n$ can fail. Inserting the upper bound on $n$ from (4.7), the estimate (4.8) follows once $N$ is sufficiently large. $\quad\square$

Our next goal is to show that the event $\mathcal{B}_n^c$ is quite unlikely to occur:

LEMMA 4.3.    *Let* $\gamma \in (s/(2d), 1)$ *and let* $s' \in (s, 2d\gamma)$. *Let* $N = |x - y| \geq N'$ *where* $N'$ *is as in Lemma* 4.2 *and define* $N_k = N^{\gamma^k}$. *Then there is a constant* $c_3 > 0$ *such that if* $\ell$ *in the definition of* $\mathcal{B}_k$ *is an odd integer between* $N_n$ *and* $2N_n$, *then*

$$(4.17) \qquad \mathbb{P}(\mathcal{B}_{k+1}^c \cap \mathcal{B}_k) \leq 2^{k+1} \exp\{-c_3 N_k^{2d\gamma - s'}\},$$

*for all* $k < n$, *where* $n$ *is as in Lemma* 4.2. *In particular,*

$$(4.18) \qquad \mathbb{P}(\mathcal{B}_n^c) \leq 2^{n+1} \exp\{-c_3 N_n^{2d\gamma - s'}\}.$$

PROOF.    Clearly, (4.18) is a result of summing (4.17), so we just need to prove (4.17) for all $k = 0, 1, \ldots, n$. By the fact that $N \geq N'$, we can assume that the scales $N_k$ and $N_{k+1}$ are clearly separated in the sense of the inequalities (4.11) and (4.13). Let $\mathcal{B}'_m$ be the event that there exists a hierarchy $\mathcal{H}_m(x, y)$ of depth $m$ connecting $x$ and $y$ such that for each $k \leq m - 2$ and each $\sigma \in \{0, 1\}^k$,

$$(4.19) \qquad z_{\sigma 01} \in B_{N_{k+1}}(z_{\sigma 0}) \quad \text{and} \quad z_{\sigma 10} \in B_{N_{k+1}}(z_{\sigma 1}).$$

A comparison with (4.5) shows that $\mathcal{B}_k \subset \mathcal{B}'_k$. Consider also the following events:

1. The event $\mathcal{A}_1$ that, for any hierarchy $\mathcal{H}_k(x, y)$ that would make $\mathcal{B}'_k$ satisfied, we have $|\mathscr{D}_{N_k}^{(\rho,\ell)}(z_\sigma)| \leq \rho N_k^d$ for some $\sigma \in \{0, 1\}^k$.

2. The event $\mathcal{A}_2$ that, for any hierarchy that would make $\mathcal{B}'_k$ satisfied, there exists a pair of sites $(z, z')$ of the type $(z_{\sigma 0}, z_{\sigma 01})$ or $(z_{\sigma 1}, z_{\sigma 10})$ with $\sigma \in \{0, 1\}^{k-2}$ such that there is no occupied bond between the sets $\mathscr{D}_{N_k}^{(\rho,\ell)}(z)$ and $\mathscr{D}_{N_k}^{(\rho,\ell)}(z')$.

Now on $\mathcal{B}_k \cap \mathcal{B}_{k+1}^c$ there exists a hierarchy that would make $\mathcal{B}'_k$ satisfied, but such that for some pair of sites as in the definition of $\mathcal{A}_2$, the sets $\mathscr{D}_{N_k}^{(\rho,\ell)}(z)$ and $\mathscr{D}_{N_k}^{(\rho,\ell)}(z')$ are not connected by an occupied bond. It follows that $\mathcal{B}_k \cap \mathcal{B}_{k+1}^c \subset \mathcal{B}'_k \cap \mathcal{A}_2$. The event $\mathcal{A}_1$ will be used to write $\mathcal{A}_2$ as the union of $\mathcal{A}_1$ and $\mathcal{A}_1^c \cap \mathcal{A}_2$, whose probabilities are more convenient to estimate.

The proof now mimics the argument from the proof of Lemma 4.2. By the fact that $N \geq N'$, the event $\mathcal{B}'_k$ is determined by looking only at the bonds that are longer than $b N_{k-2}$. Let $\mathscr{F}'$ denote the $\sigma$-algebra generated by the random variables $(\omega_{zz'})$ with $|z - z'| \geq b N_{k-2}$. Then $\mathcal{B}'_k \in \mathscr{F}'$. On the other hand, conditional on $\mathscr{F}'$, the event $\mathcal{A}_1$ is only determined by looking at the bonds that are shorter than $\operatorname{diam} \Lambda_\ell$. By (4.3), we have

$$(4.20) \qquad \mathbb{P}(\mathcal{A}_1 | \mathscr{F}') \leq 2^k \exp\{-\rho N_k^{2d-s'}\} \qquad \text{on } \mathcal{B}'_k.$$

Here $2^k$ counts the number of pairs where $\mathcal{A}_1$ can go wrong.

Concerning the event $\mathcal{A}_2$, we note that conditional on $\mathcal{A}_1^c \cap \mathcal{B}_k$, the event $\mathcal{A}_2$ is determined by the bonds of length between $b N_{k-1}$ and $b^{-1} N_{k-1}$, which by (4.13) must be either longer than $\operatorname{diam} \Lambda_\ell$ or shorter than $b N_{k-2}$. Let $\mathscr{F}$ be the $\sigma$-algebra generated by $(\omega_{zz'})$ with $|z - z'| \geq b N_{k-2}$ or $|z - z'| \leq \operatorname{diam} \Lambda_\ell$. Then $\mathcal{A}_1^c \cap \mathcal{B}'_k \in \mathscr{F}$ and $\mathcal{A}_2$ is determined by bonds independent of $\mathscr{F}$. Let $\beta$ be the same constant as in the proof of Lemma 4.2. Then we have

$$(4.21) \qquad \mathbb{P}(\mathcal{A}_2 | \mathscr{F}) \leq 2^k \exp\left\{-\beta \rho^2 \frac{N_k^{2d}}{(N_{k-1}/b)^s}\right\} \qquad \text{on } \mathcal{A}_1^c \cap \mathcal{B}_k.$$

Putting these bounds together and choosing $c_3$ appropriately, (4.17) directly follows. $\quad\square$

Lemmas 4.2 and 4.3 finally allow us to prove Proposition 4.1.

PROOF OF PROPOSITION 4.1. Let $\Delta' > \Delta$ and let $\varepsilon \in (0, 1)$. Choose $\gamma \in (s/(2d), 1)$ such that (4.6) holds true and pick an $s' \in (s, 2d\gamma)$. Suppose $N \geq N'$, where $N'$ is the constant from Lemmas 4.2 and 4.3, and let $n$ be as in Lemma 4.2. Fix an odd integer $\ell$ between $N_n$ and $2N_n$ and let $c_3$ be the constant from Lemma 4.3.

Invoking the inclusion

(4.22)
$$\{D(x,y) \geq (\log N)^{\Delta'}\}$$
$$\subset (\{D(x,y) \geq (\log N)^{\Delta'}\} \cap \mathcal{B}_n \cap \mathcal{T}) \cup \mathcal{B}_n^c \cup \mathcal{T}^c,$$

we just need to estimate the probability of the three events on the right-hand side. Lemma 4.2 shows that the probability of the first event is less than $\varepsilon$. Lemma 4.3 in conjunction with the bounds (4.9) shows that

(4.23)
$$\mathbb{P}(\mathcal{B}_n) \leq 2 \exp\{\Delta'(\log \log N) - c_3 e^{(2d\gamma - s')(\log \log N)^\varepsilon}\},$$

which can also be made less than $\varepsilon$ by choosing $N$ sufficiently large. Finally, the probability $\mathbb{P}(\mathcal{T}^c)$ is estimated from (4.4) where we assume that $N$ is so large that also $\ell \geq N_n \geq \ell_1$. From here the desired claim follows.  □

4.2. *Lower bound.* The goal of this section is to prove the lower bound in (1.6). As in Section 4.1, we formulate the relevant claim as a separate proposition:

PROPOSITION 4.4. *Suppose that* (1.3) *holds with an* $s \in (d, 2d)$ *and let* $\Delta = \Delta(s,d)$ *be as in* (1.2). *For each* $\Delta' < \Delta$ *and each* $\varepsilon > 0$, *there exists an* $N_0 < \infty$ *such that*

(4.24)
$$\mathbb{P}(D(x,y) \leq (\log |x - y|)^{\Delta'}) \leq \varepsilon$$

*holds for all* $x, y \in \mathbb{Z}^d$ *with* $|x - y| \geq N_0$.

In conjunction with Proposition 4.1, this result immediately implies Theorem 1.1.

PROOF OF THEOREM 1.1. Let $\varepsilon > 0$ and let $\mathcal{D}_\varepsilon(x)$ be the event in (1.6). Choosing $\Delta'$ such that $|\Delta - \Delta'| \leq \varepsilon$, Propositions 4.1 and 4.4 ensure that $\lim_{|x| \to \infty} \mathbb{P}(\mathcal{D}_\varepsilon(x)^c \cap \{0, x \in \mathscr{C}_\infty\}) = 0$. Then (1.6) follows by noting that, by FKG inequality and translation invariance, we have $\mathbb{P}(0, x \in \mathscr{C}_\infty) \geq \mathbb{P}(0 \in \mathscr{C}_\infty)^2$, which is positive by our assumption that there is percolation.  □

The remainder of this section will be spent on the proof of Proposition 4.4. As discussed in Section 2.2, our strategy will be to show that each path connecting $x$ and $y$ in less than $(\log |x - y|)^\Delta$ steps contains a hierarchy whose "gaps" obey the conditions (2.7). (As far as this claim is concerned, the specific choice of the exponent plays no essential role; any positive number will do.) This will be used to control the combined length of the paths needed to span the "gaps" and show that it will eventually exceed $(\log |x - y|)^{\Delta'}$ for any $\Delta' < \Delta$.

We begin by defining the relevant events. Let $x, y \in \mathbb{Z}^d$ be distinct (and distant) sites and let $N = |x - y|$. Fix a number $\gamma \in (0, s/(2d))$ and, for each integer $n \geq 2$, let $\mathcal{E}_n = \mathcal{E}_{n,\gamma}(x, y)$ be the event that *every* hierarchy $\mathcal{H}_n(x, y)$ of depth $n$ connecting $x$ and $y$ such that

$$(4.25) \qquad |z_{\sigma 01} - z_{\sigma 10}| \geq |z_{\sigma 0} - z_{\sigma 1}| (\log N)^{-\Delta}$$

holds for all $k = 0, 1, \ldots, n-2$ and all $\sigma \in \{0,1\}^k$ will also satisfy the bounds

$$(4.26) \qquad \prod_{\sigma \in \{0,1\}^k} |z_{\sigma 0} - z_{\sigma 1}| \vee 1 \geq N^{(2\gamma)^k}$$

for all $k = 1, 2, \ldots, n-1$. Here "$\vee$" is a shorthand for maximum.

REMARK 4.1.   Since we allow the possibility of "site collapse" in our definition of a hierarchy (e.g., we do not forbid that $z_{\sigma 00} = z_{\sigma 01}$), we must use a "$\vee$" on the left-hand side of (4.26). Note that (4.26) is a precise from of (2.7) while (4.25) is a precise from of (2.6).

Our first goal is to estimate the probability of $\mathcal{E}_n^c$:

LEMMA 4.5.   *Let $\gamma \in (0, s/(2d))$ and let $s' \in (2d\gamma, s)$ be such that $s' > d$. Let $\mathcal{E}_n = \mathcal{E}_{n,\gamma}(x, y)$ be as above. Then there exists a constant $c_4 \in (0, \infty)$ such that for all $x, y \in \mathbb{Z}^d$ with $N = |x - y|$ satisfying $\gamma^n \log N \geq 2(s' - d)$,*

$$(4.27) \qquad \mathbb{P}(\mathcal{E}_{n+1}^c \cap \mathcal{E}_n) \leq (\log N)^{c_4 2^n} N^{-(s' - 2d\gamma)(2\gamma)^n}.$$

The proof of Lemma 4.5 requires certain combinatorial estimates whose precise statements and proofs have been deferred to Lemmas A.1 and A.2 in the Appendix. We encourage the reader to skim through the statements of these lemmas before plunging into the forthcoming proof.

PROOF OF LEMMA 4.5.   On $\mathcal{E}_{n+1}^c \cap \mathcal{E}_n$, there exists a hierarchy $\mathcal{H}_n(x, y)$ such that the bound (4.26) holds for all $k = 1, \ldots, n-1$ but does *not* hold for $k = n$. In order to estimate the probability of such an event, let $\Theta(n)$ be the collection of all $2^n$-tuples $(z_\sigma)$ of sites such that (4.25) holds for all $\sigma \in \{0,1\}^k$ with $k = 0, 1, \ldots, n-1$, while (4.26) is true only for $k = 1, \ldots, n-1$ but not for $k = n$. Then we can write

$$(4.28) \qquad \mathbb{P}(\mathcal{E}_{n+1}^c \cap \mathcal{E}_n) \leq \sum_{(z_\sigma) \in \Theta(n)} \prod_{k=0}^{n-1} \prod_{\sigma \in \{0,1\}^k} p(z_{\sigma 01}, z_{\sigma 10}),$$

where $p(z, z') = 1 - e^{-q(z-z')}$ for $z \neq z'$—see (1.4)—while $p(z, z') = 1$ for $z = z'$. As specified in the definition of the hierarchy, none of the bonds $(z_{\sigma 01}, z_{\sigma 10})$ may appear more than once, whence (4.28) follows by invoking inclusion–exclusion and independence.

In order to estimate the right-hand side of (4.28), we will introduce a convenient change of variables: For each $k = 0, 1, \ldots, n$ and each $\sigma \in \{0,1\}^k$, let

$$(4.29) \qquad\qquad t_\sigma = z_{\sigma 0} - z_{\sigma 1}.$$

(Thus, $t_\varnothing$ is just $x - y$, while $t_0$ represents the "gap" $z_{00} - z_{01}$ and $t_1$ represents the "gap" $z_{10} - z_{11}$. Note that the $N_\sigma$'s from Section 2.2 are related to the $t_\sigma$'s via $N_\sigma = |t_\sigma|$.) Clearly, once $x$ and $y$ are fixed and $t_\sigma$ are defined for all $\sigma \in \{0,1\}^k$ and all $k = 1, \ldots, n$, all of $z_\sigma$ with $\sigma \in \{0,1\}^{n+1}$ can be reconstructed from (4.29). In terms of the $t_\sigma$'s, the conditions (4.26) can be written as

$$(4.30) \qquad\qquad \prod_{\sigma \in \{0,1\}^k} |t_\sigma| \vee 1 \geq N^{(2\gamma)^k},$$

which on $\mathcal{E}_{n+1}^c \cap \mathcal{E}_n$ is required to hold for all $k = 1, 2, \ldots, n-1$ and to fail for $k = n$, while (4.25) can be rewritten as

$$(4.31) \qquad\qquad |z_{\sigma 01} - z_{\sigma 10}| \geq |t_\sigma|(\log N)^{-\Delta},$$

which is required to hold for all $k = 0, 1, \ldots, n-1$ and all $\sigma \in \{0,1\}^k$.

The latter condition—(4.31)—allows us to recast (4.28) entirely in terms of the $t_\sigma$'s. Indeed, let $\Omega(k)$ be the set of all collections $(t_\sigma)$, $\sigma \in \{0,1\}^k$, of elements from $\mathbb{Z}^d$ such that (4.30) holds true and let $C \in (0, \infty)$ be a constant so large that

$$(4.32) \qquad\qquad p(z, z') \leq \frac{C}{|z - z'|^{s'}}$$

is true for all distinct $z, z' \in \mathbb{Z}^d$. Then (4.32) and (4.30) allow us to write

$$(4.33)$$
$$\mathbb{P}(\mathcal{E}_{n+1}^c \cap \mathcal{E}_n)$$
$$\leq \sum_{(t_\sigma) \in \Omega(1)} \cdots \sum_{(t_\sigma) \in \Omega(n-1)} \sum_{(t_\sigma) \notin \Omega(n)} \prod_{k=0}^{n-1} \prod_{\sigma \in \{0,1\}^k} \frac{C(\log N)^{s'\Delta}}{(|t_\sigma| \vee 1)^{s'}},$$

where $t_\varnothing = x - y$ and where we assumed that $C$ is so large that the last fraction exceeds 1 whenever $|t_\sigma| \leq 1$.

The right-hand side of (4.33) is now estimated as follows: First we will extract the terms $C(\log N)^{s'\Delta}$ and write the sequence of sums as a product by grouping the corresponding $t_\sigma$'s with their sum (and noting that $|t_\varnothing| = N$). This gives

$$(4.34)$$
$$\mathbb{P}(\mathcal{E}_{n+1}^c \cap \mathcal{E}_n)$$
$$\leq \frac{[C(\log N)^{s'\Delta}]^{2^n}}{N^{s'}} \left( \sum_{(t_\sigma) \notin \Omega(n)} 1 \right) \prod_{k=1}^{n-1} \left( \sum_{(t_\sigma) \in \Omega(k)} \prod_{\sigma \in \{0,1\}^k} \frac{1}{(|t_\sigma| \vee 1)^{s'}} \right).$$

Now $s' > d$, which implies that the sum in the second parentheses can be estimated using Lemma A.1. Explicitly, noting that $|\{z \in \mathbb{Z}^d : n \leq |z| \vee 1 < n+1\}| \leq c(d)n^{d-1}$ for all $n \geq 1$ and some fixed constant $c(d) < \infty$, we introduce a collection of positive integers $(n_\sigma)$ and first sum over all $(t_\sigma)$ subject to the constraint $n_\sigma \leq |t_\sigma| \vee 1 < n_\sigma + 1$. Then we are in a position to apply Lemma A.1 with $\alpha = s' - d$, $b = N^{\gamma^k}$ and $\kappa = 2^k$, which yields

$$(4.35) \qquad \sum_{(t_\sigma) \in \Omega(k)} \prod_{\sigma \in \{0,1\}^k} \frac{1}{(|t_\sigma| \vee 1)^{s'}} \leq \frac{(C' \log N)^{2^k}}{N^{(s'-d)(2\gamma)^k}},$$

for some $C' < \infty$ independent of $N$ and $n$. (Here we used that $\gamma^k \leq 1$ to bound $\log N^{\gamma^k}$ by $\log N$.) The sum in the first parentheses can be estimated in a similar fashion; the result of application of Lemma A.2 with $\alpha = d$, $b = N^{\gamma^n}$ and $\kappa = 2^n$ is

$$(4.36) \qquad \sum_{(t_\sigma) \notin \Omega(n)} 1 \leq (C'' \log N)^{2^n} N^{d(2\gamma)^n}$$

for some $C'' < \infty$. Combining these estimates with (4.34) and invoking the identity (2.9), the desired bound (4.27) is proved. $\square$

Lemma 4.5 will be used to convert the event $\{D(x,y) \leq (\log N)^{\Delta'}\}$ into a statement about the total number of bonds needed to span the "gaps" of a hierarchy identified within (one of) the shortest paths connecting $x$ and $y$. Let $\mathcal{F}_n = \mathcal{F}_n(x,y)$ be the event that, for *every* hierarchy of depth $n$ connecting $x$ and $y$ and satisfying (4.25), *every* collection of (bond) self-avoiding and mutually (bond) avoiding paths $\pi_\sigma$ with $\sigma \in \{0,1\}^{n-1}$, such that $\pi_\sigma$ connects $z_{\sigma 0}$ with $z_{\sigma 1}$ without using any bond from the hierarchy, will obey the bound

$$(4.37) \qquad \sum_{\sigma \in \{0,1\}^{n-1}} |\pi_\sigma| \geq 2^n.$$

Then we have the following claim:

Lemma 4.6. *Let* $\Delta' < \Delta$. *If* $N = |x - y|$ *is sufficiently large and*

$$(4.38) \qquad n > \frac{\Delta'}{\log 2} \log \log N,$$

*then*

$$(4.39) \qquad \{D(x,y) \leq (\log N)^{\Delta'}\} \cap \mathcal{F}_n = \varnothing.$$

Proof. We will show that on $\{D(x,y) \leq (\log N)^{\Delta'}\}$ there exists a hierarchy of a depth satisfying (4.38) such that (4.25) is true, and a collection of paths $\pi_\sigma$ "spanning" the "gaps" of this hierarchy such that (4.37) is violated. Let $\pi$ be a path saturating the distance $D(x, y)$ between $x$ and $y$. The path $\pi$ is necessarily (bond) self-avoiding. Since $|\pi| \leq (\log N)^\Delta$—by our restriction to $\{D(x, y) \leq (\log N)^{\Delta'}\}$—and since $x, y \mapsto |x - y|$ satisfies the triangle inequality, the path $\pi$ must contain a bond whose length exceeds $N/(\log N)^\Delta$. Let $z_{01}$ be the endpoint of this bond on the $x$-side and let $z_{10}$ denote the endpoint on the $y$-side of the bond. Denoting $z_0 = x$ and $z_1 = y$, we thus have

$$(4.40) \qquad |z_{01} - z_{10}| \geq |z_0 - z_1| (\log N)^{-\Delta},$$

that is, (4.25) for $\sigma = \varnothing$.

Similarly, we will identify the next level of the hierarchy. Let $\pi_0$ be the portion of $\pi$ between $z_{00} = x$ and $z_{01}$, and let $\pi_1$ be the portion of $\pi$ between $z_{11} = y$ and $z_{10}$. (Note that this agrees with our notation used in the definition of the event $\mathcal{F}_n$.) Again, we have $|\pi_0|, |\pi_1| \leq (\log N)^\Delta$ and thus both $\pi_0$ and $\pi_1$ contain at least one bond of length exceeding $|z_{00} - z_{01}|/(\log N)^\Delta$ and $|z_{10} - z_{11}|/(\log N)^\Delta$, respectively. The endpoints of this bond in $\pi_0$ identify the sites $z_{001}$ and $z_{010}$, and similarly for the endpoints of the bond in $\pi_1$. (If $\pi_0$ is empty, which can only happen if $z_{01} = z_0$, we let $z_{001} = z_{010} = z_0$, and similarly for $\pi_1$.) The very construction of these bonds implies (4.25) for $\sigma = 0, 1$. Proceeding in a similar way altogether $n$-times, we arrive at a hierarchy of depth $n$ connecting $x$ and $y$ and satisfying (4.25).

The construction implicitly defines a collection of paths $\pi_\sigma$ with $\sigma \in \{0, 1\}^{n-1}$ such that $\pi_\sigma$ is the portion of $\pi$ connecting the endpoints of the "gap" $(z_{\sigma 0}, z_{\sigma 1})$. Now we are ready to prove (4.39). Indeed, if $\{D(x, y) \leq (\log N)^{\Delta'}\}$ occurs, the combined length of all $\pi_\sigma$'s must be less than $(\log N)^{\overline{\Delta'}}$, which by (4.37) is strictly less than $2^n$. But then there exists a hierarchy of depth $n$ and self-avoiding and mutually avoiding paths $\pi_\sigma$ "spanning" its "gaps" such that (4.37) is violated. Consequently, we must have $\{D(x, y) \leq (\log N)^{\Delta'}\} \subset \mathcal{F}_n^c$.   $\square$

In light of Lemma 4.6, to prove Proposition 4.4, we will need a bound on the probability of $\mathcal{F}_n^c$ for some $n$ obeying (4.38). However, invoking also Lemma 4.5, we can as well focus just on the event $\mathcal{F}_n^c \cap \mathcal{E}_n$.

LEMMA 4.7.   *Let $\gamma \in (0, s/(2d))$ and let $s' \in (2d\gamma, s)$ be such that $s' > d$. Let $\mathcal{E}_n = \mathcal{E}_{n,\gamma}(x, y)$ and $\mathcal{F}_n = \mathcal{F}_n(x, y)$ be as above. Then there exists a constant $c_5 \in (0, \infty)$ such that for all distinct $x, y \in \mathbb{Z}^d$ with $N = |x - y|$ satisfying $\gamma^n \log N \geq 2(s' - d)$,*

$$(4.41) \qquad \mathbb{P}(\mathcal{F}_n^c \cap \mathcal{E}_n) \leq (\log N)^{c_5 2^n} N^{-s'(2\gamma)^{n-1}}.$$

Proof.   The proof will closely follow that of Lemma 4.5 so we will stay rather brief. In fact, the only essential difference is that, instead of (4.36)—which we cannot use because we are no longer in the complement of $\mathcal{E}_{n+1}$—the necessary decay for the last sum will have to be provided on the basis of the containment in $\mathcal{F}_n^{\mathrm{c}}$.

We begin by noting that on $\mathcal{F}_n^{\mathrm{c}} \cap \mathcal{E}_n$, the following events must occur:

1. There exists a hierarchy $\mathcal{H}_n(x,y)$ such that (4.25) and (4.26) hold.
2. There exists a collection of self-avoiding and mutually avoiding paths $\pi_\sigma$, with $\sigma \in \{0,1\}^{n-1}$, such that $\pi_\sigma$ connects $z_{\sigma 0}$ with $z_{\sigma 1}$ without using any bonds from $\mathcal{H}_n(x,y)$.
3. The bound (4.37) fails.

As in Lemma 4.5, we will use brute force: First we will fix a hierarchy satisfying the desired condition and try to estimate the probability that, for some collection of nonnegative integers $(m_\sigma)$, the length of the path $\pi_\sigma$ is $m_\sigma$ for each $\sigma \in \{0,1\}^{n-1}$. The fact that the paths and the hierarchy are all disjoint then allows us to write

$$(4.42) \quad \mathbb{P}\big(\mathcal{F}_n^{\mathrm{c}} \text{ occurs and } (|\pi_\sigma|) = (m_\sigma) \,|\, \mathcal{H}_n(x,y)\big) \leq \prod_{\sigma \in \{0,1\}^{n-1}} Q_{m_\sigma}(z_{\sigma 0}, z_{\sigma 1}),$$

where

$$(4.43) \qquad Q_m(z,z') = \sum_{\substack{\pi=(z_0,\dots,z_m) \\ z_0=z,\, z_m=z'}} \prod_{i=0}^{m-1} \frac{C'}{(|z_i - z_{i+1}| \vee 1)^{s'}}.$$

Here the sum runs over self-avoiding paths $\pi$ of length $m$ and $C'$ is so large that the last quotient is an upper bound on the probability that $z_i$ and $z_{i+1}$ are connected by an occupied bond.

To estimate (4.42), we first need a bound on $Q_m(z,z')$. To that end we note that, in light of the inequality $s' > d$, there exists a constant $a \in (1,\infty)$ such that for all $x,y \in \mathbb{Z}^d$,

$$(4.44) \qquad \sum_{z \in \mathbb{Z}^d} \frac{1}{(|x-z| \vee 1)^{s'}} \frac{1}{(|y-z| \vee 1)^{s'}} \leq \frac{a}{(|x-y| \vee 1)^{s'}}.$$

From here we conclude that

$$(4.45) \qquad Q_m(z,z') \leq \frac{(C'a)^m}{(|z-z'| \vee 1)^{s'}};$$

that is, up to a multiplicative factor, $Q_m(z,z')$ acts similarly as $p(z,z')$. The paths still carry some entropy in the choice of the integers $(m_\sigma)$ which amounts to counting the number $\#(m,n)$ of ordered partitions of a nonnegative integer $m$ into $2^{n-1}$ nonnegative integers. A simple estimate shows that $\#(m,n) \leq 2^{m+2^{n-1}}$ and, noting that on $\mathcal{F}_n^{\mathrm{c}}$ we only need to consider $m < 2^n$,

$$(4.46) \qquad \mathbb{P}(\mathcal{F}_n^{\mathrm{c}} | \mathcal{H}_n(x,y)) \leq (4aC')^{2^n} \prod_{\sigma \in \{0,1\}^{n-1}} \frac{1}{(|z_{\sigma 0} - z_{\sigma 1}| \vee 1)^{s'}},$$

because $\sum_{m < 2^n} \#(m,n) \leq 4^{2^n}$.

Having dispensed with the paths $\pi_\sigma$, we now start estimating the probability of $\mathcal{F}_n^{\mathrm{c}} \cap \mathcal{E}_n$. Let $\Theta^\star(n)$ be the set of all collections $(z_\sigma)$, $\sigma \in \{0,1\}^n$, obeying (4.25) for $k = 0, 1, \ldots, n-2$ and (4.26) for $k = 1, \ldots, n-1$. The bounds (4.32), (4.25) and (4.46) then give

$$(4.47) \quad \mathbb{P}(\mathcal{F}_n^{\mathrm{c}} \cap \mathcal{E}_n) \leq (4aC')^{2^n} \sum_{(z_\sigma) \in \Theta^\star(n)} \prod_{k=0}^{n-1} \prod_{\sigma \in \{0,1\}^k} \frac{(C \log N)^{s' 2^n}}{(|z_{\sigma 0} - z_{\sigma 1}| \vee 1)^{s'}}.$$

Here $C$ is the same constant as in (4.33) and the product still goes only up to $(n-1)$—despite the insertion of the terms from (4.46)—because we are now looking only at a hierarchy of depth $n$ (and not $n+1$ as in the proof of Lemma 4.5). Passing again to the variables $t_\sigma = z_{\sigma 0} - z_{\sigma 1}$ and recalling the definition of $\Omega(k)$ from the proof of Lemma 4.5, we now get

$$(4.48) \quad \mathbb{P}(\mathcal{F}_n^{\mathrm{c}} \cap \mathcal{E}_n) \leq \frac{[C''(\log N)^{s'\Delta}]^{2^n}}{N^{s'}} \prod_{k=1}^{n-1} \left( \sum_{(t_\sigma) \in \Omega(k)} \prod_{\sigma \in \{0,1\}^k} \frac{1}{(|t_\sigma| \vee 1)^{s'}} \right),$$

where $C'' < \infty$. Each term in the product can now be estimated by (4.35). Using

$$(4.49) \qquad s' + (s' - d) \sum_{k=1}^{n-1} (2\gamma)^k \geq s'(2\gamma)^{n-1}$$

instead of (2.9), the estimate (4.41) directly follows.    $\square$

Having assembled all necessary ingredients, we can now finish the proof of Proposition 4.4.

PROOF OF PROPOSITION 4.4. Let $\Delta' < \Delta$ and, recalling that $2^{-1/\Delta} = s/(2d)$, choose an $s' \in (d, s)$ such that $2^{-1/\Delta'} < s'/(2d)$. Pick a number $\gamma$ such that

$$(4.50) \qquad 2^{-1/\Delta'} < \gamma < \frac{s'}{2d}$$

and let $\delta = \frac{1}{2}(s' - 2d\gamma)$.

By Lemma 4.6, we have $\{D(x,y) \leq (\log N)^{\Delta'}\} \subset \mathcal{F}_n^{\mathrm{c}}$ once $n$ satisfies the bound (4.38). On the other hand, if $n$ also obeys the bound

$$(4.51) \qquad n \log(1/\gamma) \leq \log \log N - 2 \log \log \log N,$$

which is possible for large $N$ by (4.50), then we have $\gamma^n \log N \geq (\log \log N)^2$. This shows that, for $N$ large enough, the right-hand side of the bound from Lemma 4.5 is less than $N^{-\delta(2\gamma)^n}$ and similarly for the bound in Lemma 4.7. Consequently, both bounds are summable on $n$ and, increasing $N$ if necessary, the result can be made smaller than any number initially prescribed. Hence, for any $\varepsilon > 0$ and $N$ sufficiently large, we will have

$$(4.52) \qquad \mathbb{P}(\mathcal{F}_n^c) \leq \mathbb{P}(\mathcal{E}_n^c) + \mathbb{P}(\mathcal{F}_n^c \cap \mathcal{E}_n) \leq 2\varepsilon$$

once $n$ satisfies both (4.38) and (4.51). By the inclusion $\{D(x,y) \leq (\log N)^{\Delta'}\} \subset \mathcal{F}_n^c$, this finishes the proof. $\quad\square$

## APPENDIX

Here we establish the bounds needed in the proof of Lemmas 4.5 and 4.7. To that end, let $\kappa$ be a positive integer and, for $b > 0$ real, let

$$(A.1) \qquad \Xi_\kappa(b) = \left\{(n_i) \in \mathbb{N}^\kappa \colon n_i \geq 1, \prod_{i=1}^\kappa n_i \geq b^\kappa\right\}.$$

We will also use $\Xi_\kappa^\star(b)$ to denote a (formal) complement of this set, that is, the set of all collections $(n_i) \in \mathbb{N}^\kappa$ of positive integers such that $\prod_i n_i < b^\kappa$.

LEMMA A.1.    *For each $\varepsilon > 0$ there exists a constant $g_1 = g_1(\varepsilon) < \infty$ such that*

$$(A.2) \qquad \sum_{(n_i) \in \Xi_\kappa(b)} \prod_{i=1}^\kappa \frac{1}{n_i^{1+\alpha}} \leq (g_1 \, b^{-\alpha} \log b)^\kappa$$

*is true for all $\alpha > 0$, all $b > 1$ and all positive integers $\kappa$ satisfying*

$$(A.3) \qquad \alpha - \frac{\kappa-1}{\kappa \log b} \geq \varepsilon.$$

PROOF.    As is common for this kind of estimates, we will turn the sum into an integral. With each $(n_i) \in \Xi_\kappa(b)$, we will associate a unique hypercube $\mathfrak{h}(n_i) = (n_i) + [-\frac{1}{2}, \frac{1}{2})^\kappa$ in $\mathbb{R}^\kappa$ and note that whenever $(x_i) \in \mathfrak{h}(n_i)$, we have $x_i \geq n_i - 1/2 \geq n_i/2$ and $x_i \leq n_i + \frac{1}{2} \leq 2n_i$ for all $i = 1, \ldots, \kappa$. This implies that the product on the left-hand side of (A.2) can be bounded by the product of $(x_i/2)^{-(1+\alpha)}$ and

$$(A.4) \qquad \bigcup_{(n_i) \in \Xi_\kappa(b)} \mathfrak{h}(n_i) \subset \left\{(x_i) \in \mathbb{R}^\kappa \colon 2x_i \geq 1, \prod_{i=1}^\kappa (2x_i) \geq b^\kappa\right\}.$$

Noting that the $\mathfrak{h}(n_i)$ are disjoint, we can now bound the sum in (A.2) by the integral over the set on the right-hand side. Relabeling $2x_i$ by $x_i$, we

thus get

(A.5) left-hand side of (A.2) $\leq 2^{\kappa(1+2\alpha)} \underset{\prod_{i=1}^{\kappa} x_i \geq b^{\kappa}}{\underset{x_i \geq 1, i=1,\ldots,\kappa}{\int \cdots \int}} dx_1 \cdots dx_{\kappa} \prod_{i=1}^{\kappa} \frac{1}{x_i^{1+\alpha}}.$

To evaluate the integral, we introduce the substitutions $x_i = e^{y_i}$ followed by $z_j = y_1 + \cdots + y_j$ for $j = 1, \ldots, \kappa$. Since $y_i \geq 0$, the $z_j$'s are ordered and since the integrand depends only on $z_\kappa$, the integrals over $z_1, \ldots, z_{\kappa-1}$ can readily be performed. The result is

(A.6)     right-hand side of (A.5) $= 2^{\kappa(1+2\alpha)} \int_{\kappa \log b}^{\infty} dz \, \frac{z^{\kappa-1}}{(\kappa-1)!} e^{-\alpha z},$ where we have now dropped the subscript "$\kappa$" from $z$. Now the assumption (A.3) ensures that for $z \geq \kappa \log b$, the function $z \mapsto z^{\kappa-1} e^{-\alpha z}$ is strictly decreasing and, in fact, its logarithm is concave. Applying (A.3) we easily derive that for any $z \geq \kappa \log b$,

(A.7)              $z^{\kappa-1} e^{-\alpha z} \leq (\kappa \log b)^{s-1} b^{-\alpha \kappa} e^{-\varepsilon(z-\kappa \log b)}.$

Substituting this into (A.6), the integral is now easily performed. The calculation is concluded by using Stirling's formula to cancel the factor $\kappa^{\kappa-1}$ coming from the previous estimation against the leading order of $(\kappa-1)!$ in the denominator. $\square$

Our next claim concerns a similar sum over the indices in $\Xi_\kappa^\star(a)$:

LEMMA A.2.   *There exists a constant $g_2 < \infty$ such that for each $\alpha \geq 1$, each $b \geq e/4$ and any positive integer $\kappa$,*

(A.8)                     $\sum_{(n_i) \in \Xi_\kappa^\star(b)} \prod_{i=1}^{\kappa} n_i^{\alpha-1} \leq (g_2 \, b^\alpha \log b)^\kappa.$

PROOF.   A moment's thought reveals that we only have to address the case $\alpha = 1$. We will call upon the argument from Lemma A.1. Indeed, replacing (A.4) by

(A.9)        $\bigcup_{(n_i) \in \Xi_\kappa^\star(a)} \mathfrak{h}(n_i) \subset \left\{ (x_i) \in \mathbb{R}^\kappa : 2x_i \geq 1, \prod_{i=1}^{\kappa} x_i \leq (2b)^\kappa \right\},$

we easily find out that

(A.10)            $\sum_{(n_i) \in \Xi_\kappa^\star(a)} 1 \leq 2^{-\kappa} \underset{\prod_{i=1}^{\kappa} x_i \leq (4b)^\kappa}{\underset{x_i \geq 1, i=1,\ldots,\kappa}{\int \cdots \int}} dx_1 \cdots dx_\kappa.$

Invoking the same substitutions as before, we then get that the right-hand side of (A.10) equals

$$(A.11) \qquad 2^{-\kappa} \int_0^{\kappa \log(4b)} dz\, \frac{z^{\kappa-1}}{(\kappa-1)!}\, e^z \le 2^{-\kappa} (4b)^\kappa \frac{1}{\kappa!} (\kappa \log(4b))^\kappa.$$

Here we used the bound $e^z \le (4b)^\kappa$ to get rid of the exponential in the integral and then integrated out. Invoking Stirling's formula, the desired bound directly follows. $\quad\square$

**Acknowledgments.** This work was started when the author was a member of Theory Group in Microsoft Research, Redmond, WA. Discussions with Noam Berger and Itai Benjamini are gratefully acknowledged. The author wishes to thank an anonymous referee for suggestions that led to improvements in the presentation.

Department of Mathematics
University of California
Los Angeles, California 90095-1555
USA
e-mail: biskup@math.ucla.edu
url: http://www.math.ucla.edu/~biskup/index.html